
\magnification=1200

\font\big=cmbx10 scaled\magstep 1
\font\huge=cmbx10 scaled\magstep 2

\def\um{U^\mu}
\def\i{\infty}
\def\loh{L^1_h(G)}
\def\loe{L^1_h({\bf E})}
\def\lob{L^1_h({\bf B})}
\def\lib{L^\infty({\bf B})}
\def\B{{\bf B}}
\def\BO{{\bf B}_0}
\def\rn{{\bf R}^n}
\def\R{{\bf R}}
\def\D{{\bf D}}
\def\E{{\bf E}}
\def\L{{\cal L}}
\def\o{\omega}
\def\O{\Omega}
\def\ep{\hfill $\triangleleft$}
\def\ef{\eqno{\triangleleft}}
\centerline{{\huge{Best Approximation in the Mean}}} \par
\medskip
\centerline{{\huge{by Analytic and Harmonic Functions}}} \par
\bigskip
\centerline{Dmitry Khavinson\footnote*{Partially supported by NSF Grant
DMS 97-03915}, 
John E. M\raise.45ex\hbox{c}Carthy\footnote\dag{Partially supported by NSF
Grant DMS 95-31967}, 
and Harold S. Shapiro} \par

\bigskip
\bigskip
\baselineskip = 18pt

{\bf Abstract: } We consider the problem of finding the best harmonic or
analytic approximant to a given function.
We discuss when the best approximant is unique, and
what regularity properties the best approximant inherits from the
original function. All our approximations are done in the mean with
respect to Lebesgue measure in the plane or higher dimensions.

\centerline{{\big{1.  Introduction.}}} \par

\bigskip
\indent
For $n \geq 2$, let $\B_n$ denote the unit ball in $\rn$, and for $p \geq 1$
let $L^p$
denote the Banach space of $p$-summable functions on $\B_n$. Let 
$L^p_h(\B_n)$ denote the subspace of harmonic functions on $\B_n$ that are
$p$-summable. When $n=2$, we often write $\D$ instead of $\B_2$, and we let 
$A^p$ denote the Bergman space of analytic functions in $L^p$.

Let $\omega$ be a function in $L^p$. We are interested in finding the best
approximation to $\omega$ in $A^p$ and $L^p_h(\B_n)$.
Existence of a best approximant is straighforward; this paper considers the
following two qualitative properties:

\smallskip
\item{(i)}   Uniqueness of best approximants, when $p =1$.
\smallskip
\item{(ii)}  Hereditary regularity of the best approximant $f^{\star}$
              inherited from that of $\omega$, e.g., whether continuity,
              \hbox{H{\"o}lder} continuity, real-analyticity of $\omega$ in
              the closed unit disk enforce those properties in $\omega$'s
              best approximant.
\par

\medskip
\indent
These and many other similar questions have been well-studied for the case
when the normalized area measure $dA := \displaystyle{{1}\over{\pi}}dxdy$ 
is replaced by
$d\sigma=\displaystyle{{d\theta}\over{2\pi}}$ on the unit circle
${{\bf{T}}}$
and the spaces $A^{p}$ are replaced, accordingly, by the familiar
\hbox{Hardy} spaces $H^{p}$ (cf., e.g., [Ak], [D], [Ka],
\hbox{[Kh2--6]}, [RS], [W],and references cited therein).  In that
situation, the approach based on \hbox{Hahn}--\hbox{Banach} duality and the
\hbox{F. and M. Riesz} theorem identifying the annihilator ${{\it{Ann}}}
\left( H^{p} \right) $ in $L^{q} \left( {{\bf{T}}},d\theta \right) $ as
$H^{q}_{0}= \left\{f\in{H}^{q}:f(0)=0 \right\} $,
$q=\displaystyle{{p}\over{p-1}}$ turns out to be quite successful and
answers a number of questions.
The difficulty with this approach when using area measure
is the tremendous size of the annihilator
${{\it{Ann}}} \left(A^{p} \right) $ of $A^{p}$ in $L^{q}$.  The following
result, which we shall call 
\hbox{Khavin's} lemma,  characterizes
${{\it{Ann}}} \left( A^{p} \right) $. \par

\medskip
\indent
For
$$
p:1<p<\infty, \hskip .12in q:{{1}\over{p}}+{{1}\over{q}}=1
$$
$$
Ann \left( A^{p} \right) := \left\{ g\in{L}^{q}:\int_{{\bf{D}}}fgdA=0
\hskip .12in {{\it{for\ all}}} \hskip .12in f\in{A}^{p} \right\}
$$
$$= \left\{ {{\partial{u}}\over{\partial\overline{z}}}\ ,
u\in{W}^{1,q}_{0}({{\bf{D}})} \right\} ,
\eqno{(1.1)}
$$
where $W^{1,q}_{0}$ is the \hbox{Sobolev} space of functions vanishing on
${{\bf{T}}}$ (cf. [KS], \hbox{[Sh 1]}).  (It can be defined as the closure
of compactly-supported test functions in $C^{\infty}_{0} \left( {{\bf{D}}}
\right) $ with respect to the \hbox{$L^{q}$-norm} of their gradients, or,
equivalently, the \hbox{$L^{q}$-norm} of their
$\displaystyle{{\partial}\over{\partial\overline{z}}}$ derivative.)  \par

\medskip
\indent
For $p=1$, one needs in (1.1) to take the \hbox{weak-$\ast$ closure} in
$L^{\infty}$ of $\displaystyle{{\partial{u}}\over{\partial\overline{z}}}$,
$u\in{C}^{\infty}_{0}({{\bf{D}}})$.  Since the dual of $L^{p}$, $p\geq{1}$,
is $L^{q}$, where $\displaystyle{{1}\over{p}}+\displaystyle{{1}\over{q}}=1$,
the general
\hbox{Hahn}--\hbox{Banach} duality relation for (1.1) then can be written in
the following form (cf. e.g., [Kh2--5], \hbox{[D, Ch.8]})
$$
\lambda\ \ :=\ \ \inf_{f\in{A}^{p}} \| \omega-f \|_{p} 
\eqno{(1.2)} 
$$
$$
\qquad\qquad =\ \max_{\matrix{g\in{{\it{Ann}}} \left( A^{p} \right) \hfill \cr
               \| g \|_{q} \leq{1} \hfill \cr }}
\left| \int_{{\bf{D}}}g\omega\ {d}A \right| .
\eqno{(1.3)}
$$
The maximum in the right side of (1.3) indicates that the extremal function
$g^{\star}\in{A}nn \left( A^{p} \right) $ always exists. \par

\medskip
\indent
The rest of the paper is organized as follows.  In \hbox{Section 2} we
prove the existence of best approximations and characterize them.  
These results are not
new (cf. [Kh2--6]), but we include them for the sake of completeness and to
set the stage for further discussion.  

Section 3 deals with the
problem of uniqueness of best approximations by analytic and
harmonic functions.  The interesting case here is, of course, $p=1$. 
We show that if $\omega$ is continuous, the best analytic approximant is
unique. For harmonic approximation, we can only show that in dimension $2$
two different best harmonic approximants to a continuous function on the open
disk cannot differ by a bounded function.

 In
\hbox{Section 4} we prove two results concerning hereditary smoothness of
best approximation by $A^p$ functions, and discuss some open problems.

Section 5 deals with 
``badly approximable'' functions.  In the harmonic case, this leads to
questions concerning harmonic peak sets, which we investigate.

In \hbox{Section 6}, we give a new proof of the theorem of
Armitage, Gardiner,  Haussmann and Rogge [AGHR] characterizing best
approximation in $L^1$ 
to functions continuous on $\overline{\B_n}$ and subharmonic on
$\B_n$ by functions continuous on $\overline{\B_n}$ and harmonic on
$\B_n$.

Finally in Section 7, 
we consider best approximation in $L^1$ to the Newton kernel. We give an
explicit example of a smooth function, real-analytic on $\partial \B_n$,
whose best harmonic approximant is unbounded. (The first example  of this
type was
given in [GHJ], where the authors showed that the best harmonic
approximant to the monomial $x^4y^4$ in $L^1(\D)$ is not continuous on
$\overline{\D}$).
We also construct a continuous function on the closed disk whose best
analytic approximant is unbounded.

Although we carry out the presentation for analytic functions in the unit disk
${{\bf{D}}}$, a large portion of the results readily extend to arbitrary
smoothly bounded domains in ${{\bf{C}}}$ with merely cosmetic 
changes to the proofs. Let us also point out that somewhat related
topics are discussed in a paper by Vukoti\'c [V].

The authors thank Makoto Sakai for a valuable communication that showed 
that some of our 
more optimistic conjectures were false, and Stephen Gardiner for pointing out
a gap in an earlier version of the paper.

\bigskip\bigskip

\centerline{{\big{2.  Existence of best approximations.}}} \par

\bigskip
\noindent
{{\bf{Theorem 2.1.}}}  {{\it{The extremal function $f^{\star}$ giving the
best approximation in (1.2) always exists.}}} \par

\medskip
\noindent
{{\bf{Proof.}}}  Let $p\geq{1}$ and let $ \left\{ f_{n} \right\} \in{A}^{p}$
be a minimizing sequence, i.e., $ \left\| \omega-f_{n} \right\|_{p}
\to\lambda$.  Then $ \left\| f_{n} \right\|_{p} \leq{C}<+\infty$ for all
$n$ and hence, by subharmonicity, $f_{n}$'s are uniformly bounded on
compact subsets of ${{\bf{D}}}$.  Therefore, taking a subsequence, we can
assume that $ \left\{ f_{n} \right\} $ converge uniformly on compact subsets
of ${{\bf{D}}}$ to $f^{\star}\in{A}^{p}$.  By \hbox{Fatou's} lemma
$$
\lambda^{p}\leq \left\| \omega-f^{\star} \right\|^{p}_{p}
\leq{\liminf }_{n\to\infty} \left\| \omega-f_{n} \right\|^{p}_{p}
=\lambda^{p}
$$
and, hence, $f^{\star}$ is the best approximant to $\omega$. Q.E.D. 
\ep

\medskip
\noindent
The following result is based on the \hbox{Hahn}--\hbox{Banach} theorem (cf.
\hbox{[D, Ch.8]}, \hbox{[Kh2--5]}, [RS]) and provides the standard necessary
and sufficient conditions for the functions $f^{\star},g^{\star}$ to be the
extremals in the respective problems~(1.2),~(1.3). \par

\medskip
\noindent
{{\bf{Theorem 2.2.}}}  {{\it{
\item{(i)}   Let $p>1$,
             $q:\displaystyle{{1}\over{p}}+\displaystyle{{1}\over{q}}=1$.
             $f^{\star}\in{A}^{p},g^{\star}\in{{\it{Ann}}} \left( A^{p}
             \right) $ are extremals in (1.2) and (1.3) if and only if,
             for some $\alpha\in{{\bf{R}}}$,
$$
             e^{i\alpha}g^{\star} \left( \omega-f^{\star} \right) \geq{0}
             \hskip .12in {{\it{a.e.\ in\ }}} {{\bf{D}}}
$$
             and
$$
             \lambda^{p} \left| g^{\star} \right|^{q} = \left|
             \omega-f^{\star} \right|^{p} \hskip .12in {{\it{a.e.\ in\ }}}
             {{\bf{D}}},
\eqno{(2.3)}
$$
             where $  \lambda:={{\it{dist}}}_{L^{p}} \left(
             \omega,A^{p} \right) $.
\item{(ii)}  For $p=1$, (2.3) becomes
$$
             e^{i\alpha}g^{\star} \left( \omega-f^{\star} \right) = \left|
             \omega-f^{\star} \right| \hskip .12in {{\it{a.e.\ in\ }}}
             {{\bf{D}}},
\eqno{(2.4)}
$$
             where $f^{\star}\in{A}^{1},g^{\star}\in{{\it{Ann}}} \left(
             A^{1} \right) $.
}}} \par

\medskip
\indent
For the reader's convenience we shall sketch a (standard) proof of
\hbox{(2.3)--(2.4)} (cf. \hbox{[D, Ch.8]}, \hbox{[Kh2--5]}).
\smallskip
\item{(i)}   By \hbox{Theorem 2.1} and the \hbox{Banach}--\hbox{Alaoglu}
             Theorem, there exist extremals $f^{\star},g^{\star}$.  We find,
             applying \hbox{H{\"{o}}lder's} inequality,
$$
             \lambda= \left| \int_{{\bf{D}}}g^{\star} \left(
             \omega-f^{\star} \right) dA \right| \leq\int_{{\bf{D}}} \left|
             g^{\star} \right| \left| \omega-f^{\star} \right| dA
$$
$$
             \leq \left\| g^{\star} \right\|_{q} \left\| \omega-f^{\star}
             \right\|_{p} \leq \left\| \omega-f^{\star} \right\|_{p}
             =\lambda.
\eqno{(2.5)}
$$
             Thus, equalities must occur at each step in (2.5).  Combining
             this with necessary and sufficient conditions for equality in
             \hbox{H{\"{o}}lder's} inequality we complete the proof of
             (2.3).
\smallskip
\item{(ii)}  For $p=1$, let $f^{\star}\in{A}^{1}$,
             $g^{\star}\in{L}^{\infty}\cap{{\it{Ann}}} \left( A^{1}
             \right) $ be the extremals.  Then, the chain (2.5) becomes
$$
             \lambda= \left| \int_{{\bf{D}}}g^{\star} \left(
             \omega-f^{\star} \right) dA \right| \leq\int_{{\bf{D}}} \left|
             g^{\star} \right| \left| \omega-f^{\star} \right| dA
$$
$$
             \leq \left\| \omega-f^{\star} \right\|_{L^{1}} =\lambda,
\eqno{(2.5')}
$$
             and (2.4) follows.  Conversely, if $f^{\star},g^{\star}$
             satisfy (2.3) (or, (2.4)) and $ \left\| g^{\star} \right\|_{q}
             \leq{1}$, we have equality everywhere in (2.5) (or, ($2.5'$)),
             and since for any $f\in{A}^{p}$, $g\in{{\it{Ann}}} \left( A^{p}
             \right) $, $ \left\| g \right\|_{p} \leq{1}$ we have
$$
             \left| \int_{{\bf{D}}}g \left( \omega-f^{\star} \right) dA
             \right| \leq\lambda\leq \left\| \omega-f \right\|_{p} ,
$$
             $f^{\star},g^{\star}$ must be extremals.
\ep
\medskip
Remark: For best harmonic approximation, exactly the same result holds with 
$A^p$ and $Ann(A^p)$ replaced by $L^p_h$ and $Ann(L^p_h)$.
\par
As an application of the theorem, consider the problem of approximating
the monomials $\omega=z^{n}\overline{z}^{m}$. \par

\medskip
\noindent
{{\bf{Proposition 2.3.}}}  {{\it{For $m>n$, $f^{\star}=0$.  When $n\geq{m}$,
$f^{\star}=cz^{n-m}$, where $c=c(n,m,p)$ is an appropriate constant.}}} \par

\medskip
\noindent
{{\bf{Proof.}}}  First consider the case $m>n$.  Note that
$g:=\displaystyle{{ \left| z^{n}\overline{z}^{m} \right|^{p}
}\over{z^{n}\overline{z}^{m}}} \in{{\it{Ann}}} \left( A^{p} \right) $.
(This is checked right away by going to polar coordinates.)  Hence, for
$$
g^{\star}= \left\{ \matrix{
{{g}\over{ \left\| g \right\|_{q} }} , \hfill & p>1 \hfill \cr
g                                    , \hfill & p=1 \hfill \cr
} \right.
$$
and $f^{\star}=0$, the conditions (2.3) (or, (2.4)) are satisfied and the
statement follows.  

For $n\geq{m}$, we first find $c:=c(n,m,p)$ so that
$g:=\displaystyle{{ \left| z^{n}\overline{z}^{m}-cz^{n-m} \right|^{p}
}\over{ z^{n-m} \left( \left| z \right|^{2m} -c \right) }}\in{{\it{Ann}}}
\left( A^{p} \right) $.  Note that (setting $\left| z \right| =:r$)
$g=\displaystyle{{r^{p(n-m)} \left| r^{2m}-c \right|^{p} }\over{z^{n-m}
\left( r^{2m}-c \right) }} $.  Hence, switching to polar coordinates and
integrating with respect to $\theta$ first, we observe that $g$ annihilates
all monomials $z^{k}$, $k\neq{n}-m$.  Choosing $c=c(n,m,p)<1$ so that
$$
\int^{1}_{0}r^{p(n-m)} \left| r^{2m}-c \right|^{p-1} sgn \left( r^{2m}-c
\right) rdr=0,
$$
we have $g\in{{\it{Ann}}} \left( A^{p} \right) $.  Then, as before, setting
$$
g^{\star}= \left\{ \matrix{
{{g}\over{ \left\| g \right\|_{q} }} , \hfill & p>1 \hfill \cr
g                                    , \hfill & p=1 \hfill \cr
} \right.
$$
and applying \hbox{(2.3)--(2.4)} we complete the proof. 
\ep

Remark: The same argument shows that the best harmonic approximant to $z^n
\overline{z}^m$ is $c z^{n-m}$ if $n \geq m$, and $c \overline{z}^{m-n}$ if
$m \geq n$. 
\bigskip\bigskip

\centerline{{\big{3.1.  Uniqueness of the best analytic approximation.}}}
\par

\bigskip
\indent
The following result is originally due to
\hbox{S.Ya. Khavinson} [Kh6], where it is a part of a much more general
framework.  However, for the reader's convenience we give a straightforward
independent proof. \par

\medskip
\noindent
{{\bf{Theorem 3.1.}}}  {{\it{For $p>1$, the best approximant $f^{\star}$
in (1.2) is always unique.  For $p=1$ and $\omega$ continuous in
${{\bf{D}}}$, the best approximant $f^{\star}$ is unique.  For
discontinuous $\omega$ the best approximation need not be unique.}}} \par

\medskip
\noindent
{{\bf{Proof.}}}  For $p>1$, uniqueness is an immediate consequence of
the strict convexity of $L^{p}$ (cf., e.g., [Ak]).
Let $p=1$ and $\omega$
be continuous in ${{\bf{D}}}$.  Let $f_{1},f_{2}$ be two best approximants
to $\omega$, let $g^{\star}\in{L}^{\infty}$ be the extremal in the dual
problem so the relations (2.4):
$$
\matrix{
g^{\star} \left( \omega-f_{1} \right) \hfill & = & \left| \omega-f_{1}
\right| ; \hfill \cr
                                             &   & \hfill \cr
g^{\star} \left( \omega-f_{2} \right) \hfill & = & \left| \omega-f_{2}
\right| \hfill \cr
}
\eqno{(3.2)}
$$
hold almost everywhere in ${{\bf{D}}}$.  Let us separate the following
assertions. \par

\medskip
\noindent
{{\bf{Assertion 1.}}}  {{\it{For $z\in{{\bf{D}}}$, if $ \left|
\omega(z)-f_{1}(z) \right| = \left| \omega(z)-f_{2}(z) \right| $, then
$f_{1}(z)=f_{2}(z)$.}}} \par

\medskip
\noindent
{{\bf{Proof of Assertion 1.}}}  If $\omega(z)-f_{1}(z)=0$, the conclusion is
obvious since then $\omega(z)-f_{2}(z)=0$ also.  So, suppose
$\omega(z)-f_{1}(z)\neq{0}$.  Then (since $\omega$ is assumed to be
continuous in ${{\bf{D}}}$) there is a disk $\Delta:=\Delta(z,\rho)$
centered at $z$ such that $ \left| \omega-f_{1} \right| $ and, consequently,
also $ \left| \omega-f_{2} \right| $ are positive in $\Delta$, and by (3.2)
$ \left| g^{\star} \right| \neq{0}$ almost everywhere in $\Delta$.  Thus,
(3.2) yields
$$
{{\omega-f_{1}}\over{\omega-f_{2}}}={{ \left| \omega-f_{1} \right| }\over{
\left| \omega-f_{2} \right| }} \hskip.12in {{\it{a.e.\ in\ }}}\Delta
\eqno{(3.3)}
$$
and, hence (since both sides are continuous), (3.3) holds pointwise in
$\Delta$.  In particular, (3.3) holds at $z$ and the assertion is proved.
\par

\medskip
\noindent
{{\bf{Assertion 2.}}}  {{\it{If $p(z)$ is a real-valued continuous integrable
function
in ${{\bf{D}}}$ such that $\displaystyle{\int_{{\bf{D}}}}pdA=0$, then there
is a nontrivial continuum $K\subset{{\bf{D}}}$ on which $p=0$.}}} \par

\medskip
\noindent
{{\bf{Proof of Assertion 2}}} (obvious){{\bf{.}}}  If $p\equiv{0}$
in ${{\bf{D}}}$, there is nothing to prove.  If $U_{1}:= \left\{ z:p(z)>0
\right\} \neq\emptyset$,
then $U_{2}:={{\bf{D}}}\backslash\overline{U_1}\neq\emptyset$ by the hypothesis.
Hence $\partial{U}_{1}\cap{{\bf{D}}}$ is a continuum (since it separates
points in $U_{1}$ from those in $U_{2}$) on which $p=0$. \par

\medskip
\noindent
{{\bf{End of the proof of the theorem.}}}  Let $p(z):= \left|
\omega(z)-f_{1}(z) \right| - \left| \omega(z)-f_{2}(z) \right| $.  Since
$f_{1},f_{2}$ are both best approximants to $\omega$, the hypothesis of
\hbox{Assertion 2} is satisfied.  Hence, there is a continuum $K$ on which,
according to \hbox{Assertion 1}, $f_{1}=f_{2}$ and, accordingly,
$f_{1}=f_{2}$ everywhere in ${{\bf{D}}}$. 

\medskip
\indent
The following example shows that for discontinuous $\omega$ best
approximations need not be unique (for $p=1$, of course). \par

\medskip
\noindent
{{\bf{Example 3.2.}}}  Let ${{\bf{D}}}_{0}= \left\{ z:|z|<
\displaystyle{{1}\over{\sqrt{2}}} \right\} $, so ${{\it{Area}}} \left(
{{\bf{D}}}_{0} \right) =\displaystyle{{1}\over{2}}{{\it{Area}}} \left(
{{\bf{D}}} \right) $.  Take $\omega=\chi_{{{\bf{D}}}_{0}}$, the
characteristic function of ${{\bf{D}}}_{0}$.  Then, for any
$c:0\leq{c}\leq{1}$, $f^{\star}\equiv{c}$ gives the best approximation in
$A^{1}$ to $\omega$.  Indeed, let $g^{\star}= \left\{ \displaystyle{\matrix{
1, \hfill & z\in{{\bf{D}}}_{0} \hfill \cr
-1, \hfill & z\in{{\bf{D}}}\backslash{{\bf{D}}}_{0} \hfill \cr }} \right. $.
\par

\medskip
\noindent
Obviously, $\displaystyle{\int_{{\bf{D}}}}g^{\star}z^{n}dA=0$,
$n=1,2,\ldots$, and $\displaystyle{\int_{{\bf{D}}}}g^{\star}dA=0$, as well.
Thus, $g^{\star}\in{{\it{Ann}}} \left( A^{1} \right) $ and for any
$c:0\leq{c}\leq{1}$ we have
$$
g^{\star} \left( \omega-c \right) = \left| \omega-c \right| \hskip .12in
{{\it{a.e.\ in\ }}}{{\bf{D}}}.
$$
Thus, (2.4) is fulfilled and $f^{\star}=c$ is the best approximant to
$\omega$.  The proof of the theorem is now complete. 
\ep

\medskip
\noindent
{{\bf{Remarks.}}}
\item{(i)}   The proof given of \hbox{Thm. 3.1} extends word-for-word to
              arbitrary domains, in particular, to multiply-connected
              domains.  This is in contrast with the situation in the
              \hbox{Hardy} space setting where an $H^{1}$-best
              approximation even to a real-analytic function on the boundary
              of a finitely-connected domain (in $L^{1} \left( \left| d\xi
              \right| \right)\ $ norm) need not be unique (cf.
              \hbox{[Kh2, Section 3]}).
\smallskip
\item{(ii)}  The argument given for the proof of \hbox{Thm. 3.1} is
              sufficiently flexible and, as is easily seen, extends, e.g.,
              to $\omega$, whose set of discontinuity has measure zero in
              ${{\bf{D}}}$, is relatively closed (in ${{\bf{D}}}$) and does
              not separate ${{\bf{D}}}$ \hbox{(cf. [Kh6])}.  A more general
              result of \hbox{S.Ya. Khavinson} ([Kh6]) allows to extend
              \hbox{Thm. 3.1} to functions $\omega$ with special kinds of
              discontinuities: the limit set of $\omega$ at each point of
              discontinuity is either a segment, or contains three
              noncollinear points.  Yet, the crux of all these proofs lies
              in the fact that the zero sets of analytic functions DO NOT
              separate the disk.  This raises a very intriguing question of
              finding a plausible analogue of \hbox{Thm. 3.1} for best
              harmonic approximation in the
              \hbox{$L^{1} \left( {{\bf{D}}} \right) $-metric},
              which is addressed in the following subsection.
\smallskip
\item{(iii)}  It can also be shown that if $\omega$ is quasi-continuous, as
is the case for example for functions in the Sobolev space $W^{1,1}({\bf D})$,
then the best approximant is unique. For indeed, the boundary of the set 
$$ P \ := \ 
\{z: \left|
\omega(z)-f_{1}(z) \right| > \left| \omega(z)-f_{2}(z) \right| \}
$$
is contained in 
$$
Z \ := \
\{z: \left|
\omega(z)-f_{1}(z) \right| = \left| \omega(z)-f_{2}(z) \right| \}
\ \cap \ \{z: \omega\ {\rm continuous\ at\ } z \}
$$
union a set of arbitrarily small $1$-capacity. As $\partial P$ has positive
length and therefore positive $1$-capacity, so must $Z$; so by the proof of 
Assertion 1, $f_1$ must equal $f_2$ on a set of uniqueness for analytic
functions. (For definitions of quasi-continuous and $1$-capacity, see [EG]).

\smallskip
\item{(iv)}   It follows easily from (2.3) that for $p>1$ the extremal
              function $g^{\star}$ in the dual problem (1.3) is unique
              (up to a unimodular constant factor, of course), similarly to
              the uniqueness of $f^{\star}$.  Also, (2.4) implies
              uniqueness of $g^{\star}$ (up to a unimodular constant factor
              again) provided that $\omega$ does not coincide with an
              analytic function on a set of positive measure.  The following
              example shows that if this condition is violated,
              $g^{\star}$ may not be unique.  Let $\omega= \left\{
              \matrix{ 1, \hfill & z:0\leq|z|\leq
                       \displaystyle{{1}\over{\sqrt{3}}},
                       \sqrt{\displaystyle{{2}\over{3}}}\leq|z|<1 \hfill \cr
                       0, \hfill & {{\it{elsewhere}}}. \hfill \cr }
              \right. $.  Then, taking $f^{\star}=1$ we see that for
$$
              g^{\star}_{1}= \left\{
              \matrix{ 0, \hfill & |z|\leq \displaystyle{{1}\over{\sqrt{3}}}
                       \hfill \cr
                       -1, \hfill & \displaystyle{{1}\over{\sqrt{3}}} <|z|<
                       \sqrt{\displaystyle{{2}\over{3}}} \hfill \cr
                       1, \hfill & \sqrt{\displaystyle{{2}\over{3}}}
                       \leq|z|<1 \hfill \cr }
              \right. \hskip .12in {{\it{and}}} \hskip .12in g^{\star}_{2}=
              \left\{
              \matrix{ 1, \hfill & |z|\leq\displaystyle{{1}\over{\sqrt{3}}}
                       \hfill \cr
                       -1, \hfill & \displaystyle{{1}\over{\sqrt{3}}}<|z|<
                       \sqrt{\displaystyle{{2}\over{3}}} \hfill \cr
                       0, \hfill & \sqrt{\displaystyle{{2}\over{3}}}
                       \leq|z|<1 \hfill \cr }
              \right. ,
$$
              (2.4) are satisfied and hence the dual problem has
              a  non-unique extremal
              ($g^{\star}_{1},g^{\star}_{2}\in{{\it{Ann}}} \left( A^{1}
              \right) $ since they are both radial and have the mean value
              zero over the disk.)  Once again, we remark in passing that
              the extremal function $g^{\star}$ in the dual problem in the
              \hbox{Hardy} space context is always unique \hbox{[Kh2, RS]}.
\par

\bigskip\bigskip

\centerline{{\big{3.2.  Uniqueness of the best harmonic approximation.}}}
\par

\bigskip
As before, strict convexity of
$L^{p}$ yields the uniqueness of the best harmonic approximant to a
function $\omega$ in $L^{p}_{h}(\B_n)$ for $p>1$.  For $p=1$, the complete answer is
unknown.  \hbox{Example 3.2} shows that for discontinuous $\omega$ the best
harmonic approximant need not be unique.  Similarly, the example at the
end of the previous subsection shows that uniqueness in the dual extremal
problem fails if the function $\omega$ coincides on a set of positive
measure with a harmonic function.  Whether an analogue of \hbox{Thm. 3.1}
holds for harmonic functions ({\it i.e.} whether continuous functions have
unique best harmonic approximants in $L^1$)
is unknown to the best of our knowledge.  Here,
we give some rather special results,
which extend somewhat some of those in [GHJ], where $\omega$ was
assumed to be subharmonic and real-analytic. \par

\medskip
\noindent
{{\bf{Proposition 3.3.}}}  {{\it{Let $\omega(z)=\omega \left( |z| \right) $
be a complex-valued, radial function that is continuous and integrable in
${{\bf{D}}}$.  Then the best harmonic approximant to $\omega$ in
$L^{1}({{\bf{D}}})$ is unique, and is a constant giving the minimal value in
the \hbox{one-dimensional} problem of finding the infimum
$$
\inf \left\{ \int^{1}_{0} \left| \omega(r)-c \right| rdr, \hskip .12in
c\in{{\bf{C}}} \right\} .
$$
}}} \par

\medskip
\noindent
Let us first separate the following. \par

\medskip
\noindent
{{\bf{Lemma 3.4.}}} {{\it{For $f\in{L}^{1} \left( {{\bf{D}}} \right) $, let
$f^{\sharp}(r)=\displaystyle{{1}\over{2\pi}}\displaystyle{\int^{2\pi}_{0}}f
\left( re^{i\theta} \right) d\theta$, $0<r<1$ be the mean value of $f$ over
the circle of radius $r$.  Then, for any $u\in{L}^{1}_{h} \left( {{\bf{D}}}
\right) $ we have
$$
\int_{{\bf{D}}} \left| f-u \right| dA\geq{2} \int^{1}_{0} \left|
f^{\sharp}(r)-u(0) \right| rdr.
\eqno{(3.4)}
$$
}}} \par

\medskip
\noindent
{{\bf{Proof of the Lemma.}}}  Indeed,
$$
\int_{{\bf{D}}} \left| f-u \right| dA
={{1}\over{\pi}}\int^{1}_{0} \left\{ \int^{2\pi}_{0} \left| f \left(
re^{i\theta} \right) -u \left( re^{i\theta} \right) \right| d\theta
\right\} rdr
$$
$$
\geq{{1}\over{\pi}}\int^{1}_{0} \left. \left| \int^{2\pi}_{0} \left( f
\left( re^{i\theta} \right) -u \left( re^{i\theta} \right) \right) d\theta
\right| \right\} rdr
$$
$$
=2\int^{1}_{0} \left| f^{\sharp}(r)-u(0) \right| rdr.
$$
\par

\medskip
\noindent
{{\bf{Proof of the Proposition.}}}  Observe that (3.4) becomes equality when
$h=h^{\sharp}(r)$ is radial and $u$ is a constant.  Also, note that in view
of \hbox{Assertions 1 and 2} in the proof of \hbox{Thm.3.1} that extend
mutatis mutandis to the harmonic approximation setting (in fact, to any
setting where the approximating functions are continuous), it follows that
any two best approximants always coincide on a whole continuum of points.
Thus, to finish the proof it remains to show that a strict inequality holds
in (3.4) if $u$ is not constant. \par

\medskip
\noindent
{{\bf{Lemma 3.5.}}}  {{\it{Let $u$ be a (complex-valued) harmonic function
in the closed disk.  Then, for any $r\leq{1}$, we have
$$
\left| u(0) \right| \leq{{1}\over{2\pi}}\int^{2\pi}_{0} \left| u \left(
re^{i\theta} \right) \right| d\theta
\eqno{(3.5)}
$$
and equality holds if and only if $u={{\it{const}}}\ v$, where $v$ is a
non-negative harmonic function in ${{\bf{D}}}_{r}:= \left\{ z:|z|\leq{r}
\right\} $.}}} \par

\medskip
\noindent
Indeed, since $|u|$ is subharmonic, for equality to hold
in (3.5) $u$ must have a constant argument on
$r{{\bf{T}}}:=\partial{{\bf{D}}}_{r}$ and, hence, by the \hbox{Poisson}
formula, everywhere in ${{\bf{D}}}_{r}$. \par

\medskip
\indent
Now assume $\omega$ admits a non-constant best approximant $u$.  First of
all, by \hbox{Lemma 3.4} $\omega$ also admits a constant best approximant,
namely $u(0)$.  Replacing $\omega$ by $\omega-u(0)$ we reduce the problem to
the following:  $\omega$ is a radial, continuous function whose best
approximant is zero (i.e., $\omega$ is ``badly approximable''), and  $u$ is
another, non-constant best approximant to $\omega$ with $u(0)=0$.  By
\hbox{Lemmas 3.4 and 3.5}, for each $r$ between $0$ and $1$
the function
$\omega(r)-u(rz)=k(r)v(z)$, where $k(r)$ is a unimodular constant, and $v$
is a non-negative harmonic function in ${{\bf{D}}}$ which depends on $r$.  
Thus, the range of $u$
in ${{\bf{D}}}_{r}$ lies on a half-ray passing through $\omega(r)$.
Consider two cases:
\smallskip
\item{(i)}   The range of $\omega$ contains three noncollinear points.
             (Recall that the range of $\omega$ is a {{\it{continuous}}}
             curve).  Then let $0<a<b<c<1$ be three values of $r$ such that
             $A:=\omega(a)$, $B:=\omega(b)$, $C:=\omega(c)$ form a
             non-trivial triangle.  Then, the range of $u$ in
             ${{\bf{D}}}_{a}$, ${{\bf{D}}}_{b}$, ${{\bf{D}}}_{c}$ is
             contained in half-rays through $A$, $B$, and $C$.  Hence, the
             range of $u$ in the smallest circle ${{\bf{D}}}_{a}$ must lie
             in the intersection of these three rays, i.e., it is at most a
             point, so $u$ is a constant.
\smallskip
\item{(ii)}  The range of $\omega$ is contained in a line.  Translating and
             rotating we can assume without loss of generality that $\omega$
             is real-valued and, as before, that one of its possible best
             approximants is a zero function.  If $\omega\equiv{0}$, there
             is nothing to prove.  Otherwise, $\omega$ must change sign in
             ${{\bf{D}}}$ \hbox{(cf. (2.4))}, by which we mean that there
             exist $r_{0}:0<r_{0}<1$ such that $\omega \left( r_{0} \right)
             =0$ while either to the right, or to the left from $r_{0}$
             close to $r_{0}$ $\omega$ has either positive or negative sign.
             Without loss of generality, assume that for
             $r:r_{0}-\varepsilon<r<r_{0}$ for some small $\varepsilon>0$
             $\omega$ is positive.  For all such $r$ the range of $u$ in
             ${{\bf{D}}}_{r}$ is contained in a half-ray on the real axis
             with vertex at $\omega(r)$.  Moreover, since $u(0)=0$ while
             $\omega(r)>0$, it must always be the left half-ray.  Hence,
             letting $r\to{r}_{0}-0$ we obtain that $u\leq{0}$ in
             ${{\bf{D}}}_{r_{0}}$.  But $u(0)=0$, so $u\equiv{0}$ by the
             maximum principle and the proof is now complete.
\ep

\medskip
Remark: The statement and proof of Proposition~3.3 go through with no
difficulty to radial functions on $\B_n$ (with $rdr$ replaced by
$r^{n-1}dr$).
\medskip
\noindent
{{\bf{Theorem 3.6.}}}  {\it
Let $\omega$ be a real-valued continuous
function in $L^1( {\bf{D}})$.  Then $\omega$ cannot  have two best harmonic
approximants in $L^{1}$ whose difference is bounded.}

\medskip
\noindent
{\bf{Proof.}} Suppose $h_1$ and $h_2$ are best approximants of $\omega$.
Let $f := \omega - {{1}\over{2}} (h_1 + h_2)$, and 
$h := {{1}\over{2}} (h_1 - h_2)$.
Then $f$ is continuous on ${{\bf{D}}}$ and has $0$,  $h$ and $-h$
as best approximants. We wish to prove that if $h$ is bounded then it is
identically zero.

As $\int | f| dA = \int | f+h| dA = \int |f-h|dA$, we get that
$|f| \geq |h|$ almost everywhere, and so by continuity everywhere,
 on ${\bf D}$.
Let $P = \{ z \in {\bf D} : f(z) > 0 \}$ and
$N = \{ z \in {\bf D} : f(z) < 0 \}$.
Notice that ${\bf D} \setminus P \cup N $ is
contained in the zero-set of
$h$, and is therefore of zero area.

Define $s(z)$ to be the function that is $+1$ on $P$, $-1$ on $N$, and $0$
everywhere else on ${\bf C}$. 
Note that by the harmonic analogue of (2.4) (with $\omega = f$, $f^\star =
0$,
and $g^\star = s$), the function $s$ annihilates $L^1_h$.
Therefore
$$
0 \ = \ \int_{\bf D} z^n s(z) dA \ = \ \int_P z^n dA - \int_N z^n dA .
$$
But as 
$$
\delta_{n0}\ = \ \int_{\bf D} z^n dA \ = \ \int_P z^n dA + \int_N z^n dA ,
$$
we get that
$$
\int_P z^n dA\ = \ \int_N z^n dA\ = \ {{1}\over{2}} \delta_{n0} , \eqno{(3.6)}
$$
%
%
%
%
where $\delta_{n0}$ is the Kronecker symbol.

Now, some component of either $P$ or $N$ must intersect the disk of radius
$1/\sqrt 2$.
Without loss of generality, we can assume that some component $P_0$ of $P$
does.
Then the boundary of $P_0$ cannot intersect $\partial {\bf
D}$ in a set of positive measure. For indeed, the Cauchy transform of $N$,
the function
$$
u(z) \ = \ \int_{\bf N} {{1}\over{z-w}} dA(w)
$$ 
is continuous and bounded on the entire complex plane and analytic off $N$;
by (3.6), $u(z) = {{1}\over{2z}}$ on ${\bf C}\setminus\overline{\bf D}$.

Therefore $u$ is
analytic on $P_0$ and equal to ${{1}\over{2z}}$ on $\partial P_0 \cap \partial
{\bf D}$; if this latter set were of positive Lebesgue measure, then  
$u(z)$
would equal ${1}\over{2z}$ on all of $P_0$.
As $P_0$ intersects  the disk of radius $1/\sqrt{2}$, we get that the Cauchy
transform of $N$ is greater in modulus than $1/\sqrt{2}$ at some point.
But as $N$ has area $\pi/2$, this contradicts the Ahlfors-Beurling theorem
that says that the maximum value the Cauchy transform of a set of given area
can attain is attained 
when that set is a disk and the point in question is on the
boundary (for a proof see {\it e.g.} [GK]). A calculation shows that therefore the maximum value of 
the modulus of the Cauchy transform of a set of area $\pi/2$ is $1/\sqrt{2}$.

So $\partial P_0$ is contained in the zero-set of $h$ union a null set on 
$\partial {\bf D}$ (which is perforce a null-set also with respect to
harmonic measure for $P_0$, by Nevanlinna's majorization principle for
harmonic measures). As $h$ is bounded, and vanishes almost everywhere on 
$\partial P_0$, it must be identically zero on $P_0$, and hence on the whole
disk.
\ep

\medskip
\noindent
{\bf{Remarks}}
\item{(i)} We could weaken the hypotheses of the theorem to say that $\omega$
cannot have two best harmonic approximants whose difference raised to some
power $p > 1$ has a harmonic majorant, because again the vanishing of 
$h$ almost everywhere on $\partial P_0$ 
forces it to be identically zero.
\item{(ii)}   
The theorem is false on other domains. Let $G$ be a null quadrature domain,
{\it i.e.} a domain such that the integral of every 
$L^1_h(G)$ function is zero ({\it e.g.} the half-plane - see  [Sh1]). 
Then if $h$ is any
function in $L^1_h(G)$, the function $|h|$ has all the functions $\{ ch : -1
\leq c \leq 1 \}$ as best harmonic approximants; if $h$ is also bounded, the
theorem fails.
\item{(iii)} The previous example can be translated into a remark about
weighted approximation on the unit disk, via conformal maps. By considering
the conformal map from the disk to the right half plane, we get {\it e.g.}
that, with respect to the measure $\displaystyle {{1}\over{|z-1|^4}} dA(z)$
on the unit disk, the function $|(x-1)^3 - 3(x-1)y^2|$ 
has many best harmonic approximants
that are
bounded.

\item{(iv)} The problem with generalizing the proof to higher dimensions,
using the techniques developed in Sections~5.2 and 6, is that it is not
known whether a solution of Poisson's equation with bounded data
(so $C^{2-\varepsilon}$),
can vanish along with its gradient on a set of positive measure on the
sphere.
See [BW]
where a $C^1$ example of a non-zero harmonic function that vanishes along
with its gradient on a set of positive measure is constructed.
Of course, if one knew that $\partial P_0 \cap \partial \B_n$ actually
contained an open subset of $\partial \B_n$, there would be no problem.

\item{(v)}  
Finally, we mention that we have not touched here the questions
             related to the best uniform (\hbox{Chebyshev}) harmonic 
             approximation in ${{\bf{D}}}$ $(p=\infty)$.  Some
             results concerning the best \hbox{Chebyshev} harmonic
             approximation of subharmonic functions can be found in [HKL].
\par


\centerline{{\big{4.  Hereditary regularity}}} \par
\medskip
\centerline{{\big{of the best analytic approximation in the disk.}}} \par

\bigskip
\noindent
{{\bf{Theorem 4.1.}}}  {\it{Let $\omega$ belong to the \hbox{Sobolev} space
$W^{1,p} \left( {{\bf{D}}} \right) $, $p\geq{1}$.  Then  the best
approximant $f^{\star}\in{A}^{p}$ to $\omega$ is in the \hbox{Hardy} space
$H^{p}$.  }}

\medskip
\noindent
{{\bf{Proof.}}}  Note that by Remark (iii) after Theorem~(3.1), the best
approximant is unique even for $p=1$.
First, assume $\omega$ to be real-analytic in
$\overline{{\bf{D}}}$.  Let $P_{m}= \left\{ \right.$ polynomials in
$z$ of degree $\left. \leq{m} \right\} $ and let
$\lambda_{m}:=\displaystyle{\min_{f\in{P}_{m}}} \left\| \omega-f
\right\|_{p} $.  (Since $P_{m}$ is finite dimensional, the best
approximant $f^{\star}_{m}\in{P}_{m}$ always exists.)  Obviously,
$\lambda:=\displaystyle{\min_{f\in{A}^{p}}} \left\| \omega-f \right\|
=\displaystyle{\lim_{m\to\infty}}\lambda_{m}$, 
as the polynomials are dense in $A^p$.  
Fix $m$.  Since
$\omega$ is real-analytic, $\omega-f^{\star}_{m}\neq{0}$ a.e. in
${{\bf{D}}}$ and hence, according to an analogue of (1.2) and
\hbox{Theorem 2.2}, replacing the subspace $A^{p}$ by $P_{m}$ (the proof is
the same as that of \hbox{Thm.2.2}---cf., e.g., [Kh4--5]) it follows that
$$
{{1}\over{\lambda^{p-1}}}{{ \left| \omega-f^{\star}_{m} \right|^{p}
}\over{\omega-f^{\star}_{m}}}\in{{\it{Ann}}} \left( P_{m} \right)
\hskip .12in {{\it{in}}}\ L^{q}.
\eqno{(4.1)}
$$
Then, using \hbox{Stokes'} formula and (4.1) we obtain
$$
\eqalignno{
\int_{{\bf{T}}} \left| \omega-f^{\star}_{m} \right|^{p}
{{d\theta}\over{2\pi}}= &
{{i}\over{2\pi}}\int_{{\bf{T}}} \left|
\omega-f^{\star}_{m} \right|^{p} zd\overline{z}
\cr
= & \int_{{\bf{D}}}{{\partial}\over{\partial{z}}} \left( \left(
\omega-f^{\star}_{m}) \right)^{^{{p}\over{2}}} \left(
\overline{\omega}-\overline{f}^{\star}_{m} \right)^{{p}\over{2}} \right)
zdA
+\int_{{\bf{D}}} \left| \omega-f^{\star}_{m} \right|^{p}dA
\cr
= & {{p}\over{2}}\int_{{\bf{D}}} \left[ {{ \left| \omega-f^{\star}_{m}
\right|^{p} }\over{\omega-f^{\star}_{m}}} \left(
{{\partial\omega}\over{\partial{z}}}- \left( f^{\star}_{m} \right)^{'} 
\right) +{{ \left| \omega-f^{\star}_{m} \right|^{p}
}\over{\overline{\omega}-\overline{f}^{\star}_{m}}}
{{\partial\overline{\omega}}\over{\partial{z}}} \right] z dA
+\lambda^{p}_{m} & (4.2)
\cr
= & {{p}\over{2}}\int_{{\bf{D}}} \left[ {{ \left| \omega-f^{\star}_{m}
\right|^{p} }\over{\omega-f_{m}}} \left(
z{{\partial\omega}\over{\partial{z}}} \right) + {{ \left|
\omega-f^{\star}_{m} \right|^{p}
}\over{\overline{\omega}-\overline{f}^{\star}_{m}}}z\overline{ \left(
{{\partial\omega}\over{\partial\overline{z}}} \right) } \right]
dA+\lambda^{p}_{m},
\cr }
$$
since $z\left( f^{\star}_{m} \right)^{'} \in{P}_{m}$.
Thus, applying \hbox{H{\"{o}}lder's} inequality we
obtain
$$
\int_{{\bf{T}}} \left| \omega-f^{\star}_{m} \right|^{p}
{{d\theta}\over{2\pi}}
\leq\left\{ \matrix{
                    {{p}\over{2}} \left\| \omega-f^{\star}_{m}
                    \right\|^{{p}\over{q}}_{p} \left( \left\|
                    {{\partial\omega}\over{\partial{z}}} \right\|_{p} +
                    \left\| {{\partial\omega}\over{\partial\overline{z}}}
                    \right\|_{p} \right) +\lambda^{p}_{m}, \hfill & p>1
                    \hfill \cr
                    \hfill & \hfill \cr
                    {{1}\over{2}} \left( \left\|
                    {{\partial\omega}\over{\partial{z}}} \right\|_{1} +
                    \left\|
                    {{\partial\omega}\over{\partial\overline{z}_{p}}}
                    \right\|_{1} \right) +\lambda_{m}, \hfill & p=1.
                    \hfill \cr
} \right.
\eqno{(4.3)}
$$
$ \left( \left| \omega-f^{\star}_{m} \right|^{p-1} \in{L}^{q} \right. $ and
its \hbox{ $L^{q}$-norm is $ \left. \left\| \omega-f_{m}
\right\|^{{p}\over{q}}_{p} \right) $}.  Thus, invoking standard
inequalities for \hbox{Sobolev} spaces we obtain from (4.3)
$$
\int_{{\bf{T}}} \left| f^{\star}_{m} \right|^{p} d\theta\leq{C}
\lambda^{{p}\over{q}}_{m} \left\| \omega \right\|_{W^{1,p} \left( {{\bf{D}}}
\right) } \leq{C}_{1},
\eqno{(4.4)}
$$
where $C,C_{1}$ are constants.  Thus, all \hbox{$H^{p}$ norms} of
$f^{\star}_{m}$ are uniformly bounded, so taking a subsequence we can assume
$f^{\star}_{m}$ converges to some function ${f}^{\star}\in{H}^{p}$ 
on compact subsets of ${{\bf{D}}}$
and so $ \left( \omega-f^{\star}_{m} \right) $ converges to
$\omega-f^{\star}$ pointwise in ${{\bf{D}}}$.  In both cases we have
$$
\left\| \omega-f^{\star}
\right\|_{p} \leq\liminf_{{m\to\infty}} \left\| \omega-f^{\star}_{m}
\right\|_{p} =\lim_{m\to\infty}\lambda_{m}=\lambda,
$$
so $f^{\star}$ must be the best approximant 
to $\omega$ in $A^{p}$.  Since by \hbox{Fatou's}
lemma and (4.4)
$$
\left\| f^{\star} \right\|_{H^{p}} \leq\lim_{\overline{m\to\infty}}
\left\| f^{\star}_{m} \right\|_{H^{p}} \leq{{\it{const}}} \left\| \omega
\right\|_{W^{1,p} \left( {{\bf{D}}} \right) },
\eqno{(4.5)}
$$
it follows that 
$f^{\star}\in{H}^{p}$.  This proves the theorem for real-analytic $\omega$.
Since real-analytic functions (even polynomials) are dense in the
\hbox{Sobolev} spaces, a standard approximation argument and the fact that
the estimate (4.5) depends only on the \hbox{$W^{1,p}$-norm} of $\omega$
complete the proof. 
\ep

\medskip
\noindent
{{\bf{Remarks.}}}
\item{(i)}   The idea leading to the calculation in (4.2) goes back to
              \hbox{Ryabych} [R], where it is applied in the context of
              another extremal \hbox{problem---cf. [KS]}.
\smallskip
\item{(ii)}  Essentially, the same proof shows that the best
              {{\it{harmonic}}} approximant to a \hfill \break
              \hbox{$W^{1,p}$-function} $\omega$ in $L^{p}(\D)$ belongs to the
              class $h^{p}$ \hbox{(cf. [D])}, i.e. is  representable by a
              \hbox{Poisson}--\hbox{Lebesgue} integral with an
              \hbox{$L^{p}({{\bf{T}}})$-density} for $p>1$, or by a
              \hbox{Poisson}--{Stieltjes} integral for $p=1$.  The crucial
              step in the calculation similar to (4.2) is that if
              $h^{\star}_{m}$ is the harmonic polynomial approximant to
              $\omega$ of degree $\leq{m}$, then $z
              \displaystyle{{\partial{h}^{\star}_{m}}\over{\partial{z}}}$,
              $z \displaystyle{{\partial\overline{h^{\star}_{m}}}\over{
              \partial{z}}}$ are both analytic polynomials of degree
              $\leq{m}$ ($h^{\star}_{m}$,
              $\overline{h^{\star}_{m}}$ are harmonic!), while
              $\displaystyle{{ \left| \omega-h^{\star}_{m} \right|^{p}
              }\over{\omega-h^{\star}_{m}}}$ and its conjugate both
              annihilate $A^{p}$.
\smallskip
\item{(iii)} If $\omega$ in $L^1$ has $f^\star$ as its best $A^1$
approximant, and $p$ is any function that is positive a.e., then by looking
at the signum one sees in view of Theorem~2.2 that 
$$
n \ =\ p\,\omega\ +\ (1-p)\, f^\star
$$
also has $f^\star$ as its best $A^1$
approximant. Choosing $p$
to vanish smoothly at an isolated singularity of $\o$ inside $\D$, one 
can make $n$ smoother, and then apply Theorem~4.1 to $n$.
In particular, one gets that the best analytic approximant to
$\o(z) = \displaystyle {{1}\over{z-\lambda}}$ is in $H^1$ for all $\lambda$ in $\D$.
\par

\medskip
\indent
Another type of regularity can be derived using the ideas from [Sh2].
The key idea
is \hbox{Clarkson's} inequality \hbox{(cf. [HS, p.227])}.  Let
$p:1<p\leq{2}$, $\displaystyle{{1}\over{q}}+\displaystyle{{1}\over{p}}=1$.
Then, for any $F,G\in{L}^{p} \left( {{\bf{D}}} \right) $ we have
$$
\left\| {{F+G}\over{2}} \right\|^{q}_{p} + \left\| {{F-G}\over{2}}
\right\|^{q}_{p} \leq \left( {{1}\over{2}} \| F \|^{p}_{p} +{{1}\over{2}} \|
G \|^{p}_{p} \right)^{{q}\over{p}}.
\eqno{(4.6)}
$$
For $u\in{L}^{p} \left( {{\bf{D}}} \right) $ and $\alpha$ of modulus
$1$,
denote
by $R_{\alpha}u$ the operator
$$
R_{\alpha}u(z):=u(\alpha{z}).
\eqno{(4.7)}
$$
$R_{\alpha}$ is an isometry of $L^{p}$ and $R_{\alpha}A^{p}=A^{p}$.  Now,
let us measure the ``smoothness'' of the function (the ``mean
\hbox{H{\"{o}}lder} condition'') by saying that $u\in\Lambda_{\sigma}^p$,
$\sigma>0$, if for $0\leq{t}\leq\pi$,
$$
D_{t}u:= \left\| R_{e^{it}}u+R_{e^{-it}}u-2u \right\|_{p} =O \left(
t^{\sigma} \right) .
$$
(For simplicity, we shall just write $\Lambda_\sigma$ when the choice of
$p$ is understood). 
Of course, e.g., $u\in{C}^{2}\Rightarrow{u}\in\Lambda_{2}$, etc.

\medskip
\noindent
{{\bf{Theorem 4.2.}}}  {\it Let $ 1 < p \leq 2$, and let $q = p/(p-1)$.
Let
$\omega\in\Lambda_{\sigma}$ for some $\sigma>0$ and let $f^{\star}$ be
its
best approximant in $A^{p}$.
Then
$f^{\star}\in\Lambda_{{\sigma}\over{q}}$.} \par

\medskip
\noindent
{\bf Proof.} By scaling, we can take $\left\| \omega-f^{\star}
\right\|_{p} = 1$.
Define the operator $T_{t}$ by
$$
T_{t}(F):= \left( R_{e^{it}}F+R_{e^{-it}}F \right) /2.
\eqno{(4.8)}
$$
Clearly, $T_{t}$ is a contraction for all $t\in[0,\pi]$, and
$\omega\in\Lambda_{\sigma}$ means that
$$
\left\| T_{t}\omega-\omega \right\|_{p} \leq{C}t^{\sigma}.
$$
Let $T_{t}f^{\star}:=g\in{A}^{p}$.  Now we have, $ \left( \| \cdot \| = \|
\cdot \|_{p} \right) $:
$$\eqalignno{
\left\| g-\omega \right\| &= \left\| T_{t}f^{\star}-\omega \right\|
\cr
&\leq
 \left\|
T_{t}f^{\star}-T_{t}\omega \right\|
+ \left\| T_{t}\omega-\omega \right\| \cr
&\leq \left\| f^{\star}-\omega \right\|
+Ct^{\sigma} \cr
&= 1+Ct^{\sigma}. & (4.9)
\cr
}
$$
Applying (4.6) with $F=f^{\star}-\omega$ and $G=g-\omega$ we get
$$
\left\| {{f^{\star}+g}\over{2}}-\omega \right\|^{q} + \left\|
{{f^{\star}-g}\over{2}} \right\|^{q} \leq \left[ {{1}\over{2}} \left\|
f^{\star}-\omega \right\|^{p} +{{1}\over{2}} \left\| g-\omega \right\|^{p}
\right]^{{q}\over{p}} .
\eqno{(4.10)}
$$
Since $1= \left\| f^{\star}-\omega \right\| $, the first term on the left in
(4.10) is $\geq{1}$.  Using also (4.9) we obtain ($C$ denotes constants that
may change from line to line)
$$
1+ \left\| {{f^{\star}-g}\over{2}} \right\|^{q} \leq \left[
{{1}\over{2}}+{{1}\over{2}} \left( 1+Ct^{\sigma} \right)^{p}
\right]^{{q}\over{p}}
$$
$$
\leq \left( 1+Ct^{\sigma} \right)^{{q}\over{p}} \leq{1}+Ct^{\sigma},
\eqno{(4.11)}
$$
where $q\geq{p}$ and we take $t$ to be smaller than $1$.  Thus,
$$
\left\| f^{\star}-g \right\| \leq{C}t^{{\sigma}\over{q}}
\eqno{(4.12)}
$$
and recalling that $g=T_{t}f^{\star}$ we obtain the result.
\ep
\medskip
\noindent
{{\bf{Discussion of Theorem 4.2.}}}  Letting $\alpha=e^{it}$, the above
assertion becomes
$$
\left( \int_{{\bf{D}}} \left| f^{\star} \left( \alpha{z} \right) +f^{\star}
\left( \overline{\alpha}z \right) -2f^{\star}(z) \right|^{p} dA
\right)^{{1}\over{p}} \leq{C}t^{{\sigma}\over{q}}
$$
or, writing $f^\star (z)=\displaystyle{\sum^{\infty}_{0}}a_{n}z^{n}$,
$$
\int^{1}_{0}rdr\int^{2\pi}_{0} \left| \sum^{\infty}_{0}a_{n}r^{n} \left(
\alpha^{{n}\over{2}}-\alpha^{-{{n}\over{2}}} \right)^{2} e^{in\theta}
\right|^{p} d\theta\leq{C}t^{\sigma(p-1)}.
\eqno{(4.13)}
$$
Now, by the \hbox{Hausdorff}--\hbox{Young} inequality $(p\leq{2}!)$,
\hbox{cf. [D, p.83]}:
$$
\left( {{1}\over{2\pi}}\int^{2\pi}_{0} \left| \sum^{\infty}_{0}a_{n}r^{n}
\left( \alpha^{{n}\over{2}}-\alpha^{-{{n}\over{2}}} \right)^{2} e^{in\theta}
\right|^{p} d\theta \right)
$$
$$
\geq \left\{ \sum^{\infty}_{0} \left( \left| a_{n} \right| r^{n} \left|
\alpha^{{n}\over{2}}-\alpha^{-{{n}\over{2}}} \right|^{2} \right)^{q}
\right\}^{{p}\over{q}}.
$$
So, from (4.13) it follows $ \left( \alpha=e^{it} \right) $:
$$
\int^{1}_{0}rdr \left( \sum^{\infty}_{0} \left( \left| a_{n} \right|
r^{n}\sin^{2}{{nt}\over{2}} \right)^{q} \right)^{p-1}
\leq{C}t^{\sigma(p-1)}.
\eqno{(4.14)}
$$
Fix an integer $N$ and replace the integral on the left in (4.14) by that
over $ \left[ 1-\displaystyle{{1}\over{N}},1 \right] $ and
$\displaystyle{\sum^{\infty}_{0}}$ by $\displaystyle{\sum^{N}_{1}}$.  For
that range of $r$ and $n$, $r^{n}\geq{r}^{N}\geq \left(
1-\displaystyle{{1}\over{N}} \right)^{N} \geq{c}$, an absolute constant, so
(4.14) yields
$$
{{1}\over{N}} \left( \sum^{N}_{n=1} \left( \left| a_{n} \right|
\sin^{2}{{nt}\over{2}} \right)^{q} \right)^{p-1} \leq{C}t^{\sigma(p-1)},
$$
hence,
$$
\sum^{N}_{1} \left| a_{n} \right|^{q}
\sin^{2q}{{nt}\over{2}}\leq{C}N^{{1}\over{p-1}}t^{\sigma}.
\eqno{(4.15)}
$$
Now, for $\xi<\sigma+1$, multiply (4.15) by $t^{-\xi}$ and integrate from
$0$ to $1$.  We obtain
$$
\sum^{N}_{1} \left| a_{n} \right|^{q} \int^{1}_{0} \left(
\sin^{2}{{nt}\over{2}} \right)^{q} t^{-\xi}dt\leq{C}(\xi)N^{{1}\over{p-1}}.
\eqno{(4.16)}
$$
The integral in (4.16) is, changing variables by $nt=s$,
$$
\int^{n}_{0} \left( \sin^{2}{{s}\over{2}} \right)^{q} \left( {{n}\over{s}}
\right)^{\xi} {{ds}\over{n}}=n^{\xi-1}\int^{n}_{0} \left(
\sin^{2}{{s}\over{2}} \right)^{q}s^{-\xi}ds\geq{C}n^{\xi-1},
$$
so
$$
\sum^{N}_{1}n^{\xi-1} \left| a_{n} \right|^{q} \leq{C}(\xi)N^{{1}\over{p-1}}
$$
and
$$
\sum_{ \left[ {{N}\over{2}} \right] +1\leq{n}<N} \left| a_{n} \right|^{q}
\leq{C}(\xi)N^{{{1}\over{p-1}}-\xi+1}.
\eqno{(4.17)}
$$
(4.17) allows us to extract some particular regularity information about
$f^{\star}$.  For example, take $\sigma=2$.  Then, for all
${{3}\over{2}}<p\leq{2}$ it follows from (4.17), letting
$\delta = 2 - { 1 \over {p-1}}$ and $\tau = 3 - \xi$,
$$
\sum_{ \left[ {{N}\over{2}} \right] +1\leq{n}<N} \left| a_{n} \right|^{q}
\leq{c}(\tau)N^{(2-\delta)-(3-\tau)+1} \leq{C}(\tau)N^{-\eta}, \hskip .12in
\eqno{(4.18)}
$$
where $\eta = \delta - \tau$ is positive for
$\tau$ sufficiently small and positive.
Hence, for all $k\geq{1}$ we have
$$
\sum_{2^{k-1}\leq{n}<2^{k}} \left| a_{n} \right|^{q} <C(\tau)2^{-\eta k},
$$
so by H\"older's inequality
$$
\sum_{2^{k-1}\leq{n}<2^{k}} \left| a_{n} \right|^{2} <C(\tau)
[2^{-\eta k}]^{2 \over q} [2^k]^{{q-2} \over q}.
\eqno{(4.19)}
$$
The coefficient of $k$ in the exponent on the right-hand side of
(4.19) can, by suitable choice of $\tau$, be made negative for
$p > {8 \over 5}$, so we obtain

\medskip
\noindent
{{\bf{Corollary 4.3.}}}  {{\it{For $\displaystyle{{8}\over{5}}<p\leq
 2$ the best
approximant $f^{\star}$ in $A^{p}$ to a function $\omega\in\Lambda_{2}^p$
belongs to the \hbox{Hardy} space $H^{2}$.}}} \par

\medskip
\noindent
{{\bf{Remarks.}}}
\item{(i)}   It would be interesting to clarify the relationship between
             \hbox{Theorems 4.1 and 4.2}.
\smallskip
\item{(ii)}  \hbox{Theorems 4.1 and 4.2} certainly leave unanswered most of
             the natural regularity questions, such as: given $\omega$ to be
             real-analytic in $\overline{{\bf{D}}}$, does it imply that its
             best approximant in $A^{p}$ is \hbox{H{\"{o}}lder} continuous
             or merely continuous in $\overline{{\bf{D}}}$?  (Similar
             results and much more are known to hold for
             \hbox{$L^{p}$-approximation} on the circle--cf., e.g., [CJ],
             \hbox{[D, Ch.8]}, [Ka], [Kh2,3], [RS].) When $p=1$, much
	     regularity can be lost by harmonic aproximation - see Section~7. 
\smallskip
\item{(iii)} 
Here is another set of problems.
 For the sake of definiteness let us
take $p=1$.  Assume $\omega\in{C} \left( \overline{{\bf{D}}} \right) $, $ \|
\omega \|_{\infty} =1$ and that for some small $\varepsilon>0$ we can find
$f\in{A}^{1}$ so that
$$
\lambda:= \| \omega-f \|_{1} \leq\varepsilon.
\eqno{(4.20)}
$$
{{\bf{Question.}}}  What is the distance  from $\omega$ to the
unit ball in $H^{\infty}$ in the \hbox{$L^{1}$-norm}?
In other words,
estimate (in
terms of $\varepsilon$)
$$
\mu:=\inf \left\{ \| \omega-g \|_{1} :g\in{H}^{\infty}, \| g \|_{\infty}
\leq{1} \right\} .
\eqno{(4.21)}
$$
A similar problem in the context of \hbox{Hardy} spaces was discussed in a
recent paper \hbox{[KP-GS]}.  There, the authors showed that
\smallskip
(a) $\mu=O \left(
             \varepsilon\log\displaystyle{{1}\over{\varepsilon}} \right) $
             and
\smallskip
(b) The estimate in (a) cannot be improved to $O(\varepsilon)$.
\medskip
However, all the major ingredients of the arguments in \hbox{[KP-GS]} fail
miserably for 

Bergman functions. 
We think that the relationship between quantities (4.20) and 

(4.21) may be a
fruitful topic for future investigations.

\bigskip\bigskip

\centerline{{\big{5.1  $A^p$ Badly Approximable Functions }}}


\bigskip
\indent
We shall call a function $\omega\in{L}^{p} \left( {{\bf{D}}} \right) $
{{\it{badly approximable}}} with respect to $A^p$
if its best approximant in $A^{p}$ equals $0$.

\medskip
\noindent
{{\bf{Example 5.1.}}}
\item{(i)}   Let $a(r)\in{L}^{p} \left( rdr,[0,1] \right)$, $p\geq{1}$,
             $n\geq{1}$.  Then the function $\omega:=a(r)e^{-in\theta}$ is
             badly approximable by $A^{p}$.  Indeed,
$$
             {{ \left| \omega \right|^{p} }\over{\omega}}= \left| a(r)
             \right|^{p-1} {{\it{sgn}}}[a(r)]e^{in\theta}\in {{\it{Ann}}} 
             \left( A^{p} \right)
$$
             in $L^{q}$,
             $\displaystyle{{1}\over{p}}+\displaystyle{{1}\over{q}}=1$, and
             by \hbox{Theorem 2.2} the assertion follows.
\smallskip
\item{(ii)}  Let $p>2$, $N= \left\{ 0,z_{1},\ldots,z_{n} \right\} $---a
             finite set and let $G$ be a contractive zero divisor in
             $A^{p-2}$ \hbox{(cf. [DKSS])} corresponding to the zero set
             $N$.  Then, $G$ extends analytically across ${{\bf{T}}}$.  Set
             $\omega=\overline{G}$.  $\omega$ is badly approximable.
             Indeed, one of the characteristic properties of a contractive
             divisor is that the measure $ \left| G \right|^{p-2} dA$ is a
             representing measure for bounded analytic functions and since
             $G$ itself is bounded, for all $A^{1}$ functions.  Hence, $G
             \left| G \right|^{p-2}$ annihilates ${A}^{p} \left( G(0)=0!
             \right) $, and by \hbox{Theorem 2.2} $\omega:=\overline{G}$ is
             badly approximable.
\par

\medskip
\indent
On the other hand, we have \par

\medskip
\noindent
{{\bf{Propsition 5.2.}}}  {{\it{Let $f(z)$ be analytic and
satisfy $ \left| f(z) \right|
\geq{c}>0$ in ${{\bf{D}}}$.  Then, $\omega:=\overline{f}$ is not badly
approximable in $L^{p}$, $p\geq{1}$.}}} \par

\medskip
\noindent
{{\bf{Proof.}}}  Indeed, otherwise by \hbox{Theorem 2.2} we would have
$\displaystyle{{ \left| f \right|^{p} }\over{\overline{f}}}=f \left| f
\right|^{p-2} \bot{A}^{p}$, so in particular (as $ 
\displaystyle{{1}\over{f}}\in{H}^{\infty} $):
$$
\int_{{\bf{D}}}{{1}\over{f}}f \left| f \right|^{p-2} dA=0,
$$
an obvious contradiction. 
\ep

\medskip
\noindent
It is quite easy, using duality, to characterize all badly approximable
functions on the circle (in the context of \hbox{Hardy} spaces).  In
particular, conjugates of all inner functions vanishing at the origin are
badly approximable.  In the Bergman space, the situation is more complicated. It
can be shown, {\it e.g.}, that the functions
$$
\omega(z) \ =\ {{\overline{z}^2 (\overline{z}-a)^2}\over{(1-a \overline{z})^2}},
\quad 0< a < 1,
$$
are not badly approximable in $L^1$.
Contrast the following result with Proposition 5.2.
\medskip
\noindent
{{\bf{Proposition 5.3.}}}  {{\it{The function $\omega:= \left(
\overline{z}-a \right)^{4} $, $0<a<1$, is badly approximable in $L^{1}$.}}}
\par

\medskip
\noindent
{{\bf{Proof.}}}  In view of Theorem 2.2 and Khavin's Lemma,
 we want to find a continuous function $v$ in
$\overline{{\bf{D}}}$, $ \left. v \right|_{{\bf{T}}} =0$ such that
$$
{{\partial{v}}\over{\partial\overline{z}}}= \left[ {{(z-a)}\over{ \left(
\overline{z}-a \right) }} \right]^{2} \hskip .12in {{\it{in}}}\ {{\bf{D}}}.
\eqno{(5.1)}
$$
Integrating (5.1) we see that it is equivalent to the existence of a
holomorphic function $h$ in ${{\bf{D}}}$ satisfying
$$
h(z)=v(z)+{{(z-a)^{2}}\over{\overline{z}-a}}.
\eqno{(5.2)}
$$
On ${{\bf{T}}}$ $v=0$, so (5.2) yields that
$h(z)=\displaystyle{{z(z-a)^{2}}\over{1-az}}$, and so,
$$
v(z)={{z(z-a)^{2}}\over{1-az}}-{{(z-a)^{2}}\over{\overline{z}-a}}
$$
has all the desired properties. 
\ep

\bigskip\bigskip

\centerline{{\big{5.2  $\loh$ Badly Approximable Functions and Harmonic Peak
Sets}}}

\bigskip

Many of the ideas in harmonic approximation extend to ${\bf
R}^n$, so
we shall work there.
Let $G$ be a domain in ${\bf R}^n$, and let us introduce the following two
ways of measuring the size of a subset $F$ with respect to harmonic
functions.

{\bf Definition:} {\it For $F$ a subset of $G$, define
$$
A(F) \ = \ \sup\{ {{\int_F |h|} \over {\int_{G \setminus F} |h|}} \ : \
h \in L^1_h(G) \}
$$
and
$$
B(F) \ = \ \sup\{ {{| \int_F h|} \over {\int_{G \setminus F} |h|}} \ : \
h \in L^1_h(G) \}.
$$
}

\medskip
\noindent
{\bf Theorem 5.4.}
{\it
If $F \subseteq G$ has $B(F) > 1$, and $\omega$ in $L^1(G)$ is strictly
positive a.e. on $F$, then $\omega$ is not badly approximable.
}

\medskip
\noindent
{\bf Proof.} 
By the harmonic analogue of Theorem 2.2, if $\omega$ were badly
approximable, then there would be a function $g$ in $L^\i$ of norm one 
({\it viz.} $sgn(\overline{\omega})$) that
annihilated $\loh$ and equalled $1$ on $F$.
As $B(F) > 1$, there exists $h$ in $\loh$ such that
$|\int_F h | >  \int_{G \setminus F} | h |$.
We have
$$
\eqalignno{
\left| \int_G h \right| &= \left| \int_G h (1-g) \right| \cr
&= \left|\int_{G \setminus F} h -hg \right| \cr
&\leq 
2\int_{G \setminus F} |h| .
\cr}
$$
Now let $\lambda > 0$. Then
$$
\eqalignno{
\left|\int_G  h \right| &= \left| \int_G h (1+ \lambda g) \right| \cr
&\geq (\lambda + 1) \left|\int_F h \right| - (\lambda + 1) \int_{G \setminus F} | h |
\cr }
$$
so
$$ 
(\lambda + 1) \left|\int_F h \right| \leq
(\lambda +3) \int_{G \setminus F} | h | .
$$
Letting $\lambda \to \i$ gives
$$
\left|\int_F h \right| \leq  \int_{G \setminus F} | h |,
$$
a contradiction.
\ep

\bigskip
Note that Example 3.2 shows, with $F = {\bf D}_0$, that $B(F) = 1$ is not a
sufficient hypothesis.
 
A similar argument to the proof of Theorem 5.4 yields the following theorem.

\noindent
{\bf Theorem 5.5.} {\it
Suppose $F \subseteq G$ has $A(F) > 1$, so there is some
function $h$ in
$\loh$ such that $\int_F |h| > \int_{G \setminus F} |h|$. If $\omega$ in
$L^1(G)$ has the property that $\omega \overline{h}$ is strictly
positive a.e. on $F$, then $\omega$ is not badly approximable.
}

\bigskip
We shall call $F$ a {\it weak peak set} if $A(F) = \i$,
and a {\it strong peak set} if $B(F) = \i$.
These sets seem of interest in their own right.
A duality argument shows their connection with badly approximable functions
and dual interpolation problems.

\medskip
\noindent
{\bf Proposition 5.6.}
{\it
The set $F$ is a {weak peak set} for $L^1_h(G)$
if and only if there is a function $g$ in $L^\i(F)$ that cannot be
extended  to a bounded function on $G$  that annihilates $L^1_h(G)$.
The set $F$ is a {strong peak set} for $L^1_h(G)$
if and only if the function that is identically $1$ on $F$ cannot be
extended to a bounded function on $G$  that annihilates $L^1_h(G)$.
}

\medskip
\noindent
{\bf Proof.}
$F$ fails to be a {weak peak set} for $L^1_h(G)$
if and only if there is a constant $M$ such that
$$
\int_F |h| dA \ \leq \ M \int_{G\setminus F} |h| dA
$$
for all $h$ in $L^1_h(G)$. This implies that if $g$ is any function
in  $L^\i(F)$, then there is a function $\omega_g$ in
$L^\i(G\setminus F)$ of norm at most $M\| g \|$ such that
$$
\int_F h g dA = \int_{G\setminus F} h \omega_g dA.
$$
But then $g \chi_F - \omega_g \chi_{G\setminus F}$ is an extension of $g$ that
annihilates $L^1_h(G)$. As the reasoning is reversible, this proves the
characterization of weak peak sets.

Similarly, $F$ fails to be a {strong peak set} for $L^1_h(G)$
if and only if there is a constant $M$ such that
$$
|\int_F h dA| \ \leq \ M \int_{G\setminus F} |h| dA
$$
for all $h$ in $L^1_h(G)$.
But this implies that there is a function $\omega$ of norm at most $M$
so that
$$
\int_F h  dA = \int_{G\setminus F} h \omega dA,
$$
and so $\chi_F - \omega \chi_{G\setminus F}$ annihilates $L^1_h(G)$.
Again the argument is reversible.
\ep

\bigskip
%
%

Now we turn to geometric characterizations of peak sets, motivated by 
the previous results and 
Theorem~3.6.
\medskip
\noindent
{\bf Theorem 5.7.}
{\it
Suppose $G$ is a bounded domain in ${\bf R}^n$ and
the boundary of $G$ contains an isolated  $(n-1)$-dimensional manifold $J$
which is also in the boundary of ${\bf R}^n\setminus \overline{G}$. Then
every full neighborhood of a point in $J$ is a weak peak set for $\loh$.
}

\medskip
\noindent
{\bf Proof.}
It is easily shown that there is a point $y$ in $G$ such that a closest 
point in $\partial G$ to $y$
lies in $J$. Let $z$ be a point in $\partial G$ that is closest to $y$.
Note that the ball centered at $y$ of radius $|y-z|$ is contained in $G$.

Let $N$ be the intersection of an
open set in ${\bf R}^n$ containing $z$ with $G$.
Let $u$ be a harmonic function on ${\bf R}^n \setminus \{0\}$ with a
non-integrable singularity at $0$, such that $u$ is not integrable over any
ball with $0$ in the boundary ({\it e.g.} let $u$ be an appropriate partial
derivative of the Newton kernel).
Let $z_j$ be a sequence in ${\bf R}^n\setminus \overline{G}$ that converges
to $z$. Let $u_j(x) := u(x - z_j)$.
Then $\int_N |u_j| $ tends to infinity, while $\int_{G\setminus N} |u_j|$
stays bounded.
\ep

\bigskip
\noindent
{\bf Lemma 5.8.} {\it
Let  $G$ be a bounded domain in ${\bf R}^n$, and suppose $F$
is a weak peak set for $\loh$. Then for all $c > 0$, 
$F_c := F \cap \{x \in G : dist(x,\partial G) < c \}$ is also a weak peak
set.
}
\medskip
\noindent
{\bf Proof.} In view of the proof of Theorem~5.7, we can assume that $G
\setminus F$ has a subset $E$ of positive measure and with $cl(E) \subseteq
G$. 

As $F$ is weak peak, there is a sequence $h_j$ in $\loh$, each function
having norm one, and $\int_F |h_j| $ tending to $1$ as $j \to \i$.
Passing to a subsequence if necessary, we can assume that $h_j$ converges
uniformly on compact subsets of $G$ to a harmonic function $h$. 
As $\int_E |h_j| \to 0$, it follows that $h = 0$ on $E$ and therefore on all
of $G$. Therefore $h_j$ tends to zero uniformly on compact subsets of $G$,
and in particular on $F \setminus F_c$. So 
$\int_{F_c} |h_j| \to 1$, as desired.
\ep

\bigskip
It is possible for a set $F$ to touch the boundary but not be a weak peak
set, provided it is very thin near the boundary.

\medskip
\noindent
{\bf Theorem 5.9.} {\it
Suppose $G$ is a bounded domain in ${\bf R}^n$, and
$F \subseteq G$ satisfies 
$$
\int_F {{1}\over{(dist(z,\partial G))^n}} dA < \i .
$$
Then $F$ is not a weak peak set for $L^1_h(G)$.
}

\medskip
\noindent
{\bf Proof.}
Let $c_n$ be the volume of the unit ball in ${\bf R}^n$. For some $c > 0$,
the set $F_c$ satisfies
$$
\int_{F_c} {{1}\over{(dist(z,\partial G))^n}} dA < {{c_n}\over{2}}.
$$
By Lemma 5.8, it is sufficient to prove that $F_c$ is not a weak peak
set.
Now suppose $h$ is in $\loh$. Then by the mean value property for harmonic
functions,
$$
|h(z) | \ \leq \ {{1}\over{c_n (dist(z,\partial G)^n}} \int_G |h|dA
$$
for all $z$ in $G$. Therefore
$$
\int_{F_c} |h| dA \ \leq\ {{1}\over{2}} \int_G |h|dA
$$
so 
$$
\int_{F_c} |h| dA \ \leq\ \int_{G\setminus {F_c}} |h|dA.
\ef
$$

\bigskip
\bigskip
Characterizing strong harmonic peak sets is more subtle.
To determine whether a subset of the ball is a strong harmonic peak set, the
center is of crucial importance. Let $\bf B$ denote the unit ball in ${\bf
R}^n$, and recall that $c_n$ is its volume.
\medskip
\noindent
{\bf Theorem 5.10.}
{\it
Let $F \subseteq \B$.

(i) If $0$ is not in $\overline{F}$, then $F$ is not a strong peak
set for $\lob$.

(ii) If $0$ is in $\overline{F}$, $F$ is open and connected, and in addition
$\partial F$ contains a relatively open subset of $\partial \B$,
then $F$ is a strong peak set for $\lob$.
}

\medskip
\noindent
{\bf Proof.}
Suppose first that  $F$ omits $\B(0,r)$, the ball centered at zero of radius
$r > 0$.
Then for any integrable harmonic function $h$
$$
\int_F h  \ = c_n h(0) - \int_{\B \setminus F} h  .
$$
Therefore
$$
\eqalignno{
\left|\int_F h  \right| \ &\leq \ {{1}\over{c_n r^n}} 
\left|\int_{\B(0,r)} h  \right| + \left|\int_{\B
\setminus F} h  \right| \cr
&\leq \ {{c_n r^n + 1}\over{c_n r^n}} \int_{\B \setminus F} |h | ,
\cr }
$$
so $F$ cannot be a strong harmonic peak set.

Conversely, if $F$ is a domain that contains an open subset
 $J$ of the unit sphere in
its boundary, and
if $F$ is not a strong harmonic peak set,
let $\psi$ be a function in $L^\infty(\B)$
that annihilates $\lob$ and equals 
$1$ on $F$.

Claim: There is a $C^1$ function $u$ on ${\bf R}^n$ satisfying
$$
\Delta u = \psi, \quad u = 0 = {{\partial u}\over{\partial n}}\
{\rm on}\ J .
\eqno{(5.3)}
$$

Proof of claim: Let $E$ be the fundamental solution of the Laplacian in ${\bf
R}^n$,
and define $u$ by
$$
u\ =\ E\ast\psi .
$$
Then $\Delta u = \psi$ and $u$ is $C^1$ by elliptic regularity
[GT].
Moreover, because for $\xi \in \rn \setminus
\overline{\B}$ the function $z \mapsto E(\xi-z)$ is harmonic on $\B$,
it follows from the fact that $\psi$ annihilates $\lob$
 that $u \equiv 0$ off $\overline \B$.
As $u$ is $C^1$, 
it follows that $u$ and its first order partials vanish on
$J$.

Let $v$ be the modified Schwarz potential of $\partial B$,
{\it i.e.} the function satisfying
$$
\Delta v = 1, \quad v = 0 = {{\partial v}\over{\partial n}}\
{\rm on}\ \partial \B .
\eqno{(5.4)}
$$
As $u$ and $v$ agree on $F$ and vanish along with their gradients on $J$,
we must have $u \equiv v$ in $F$.
By direct calculation (or see [Kh1] or [Sh1]), for $n=2$ we have
$$
v(z) \ = {{1}\over{4}}(|z|^2 - 1) - {{1}\over{2}}\log|z|
$$
and for $n \geq 3$ we have
$$
v(z) \ =  {{1}\over{2n}} |z|^2 + {{1}\over{n(n-2)}} {{1}\over{|z|^{n-2}}}
- {{1}\over{2(n-2)}}.
$$ 
As $v$ has a non-removable singularity at $0$ and $u$ is bounded, 
$0$ cannot be in $\overline{F}$.
\ep
\bigskip

For $n =2$, it suffices in {\it (ii)} for $\partial F \cap \partial \D$ to
have positive measure - cf. Remark~(iv) after Theorem~3.6.

For an ellipse, the crucial points are the foci. A domain has to join
only one of these to an arc on the boundary in order to be a strong harmonic
peak set.

\bigskip
\noindent
{\bf Theorem 5.11.}
{\it
Let $\E$ be an ellipse with foci $\pm 1$, and let $F \subset \E$.

(i) If there exists a connected open set $U$ containing both foci that is
disjoint form $F$, then $F$ is not a strong peak set for $\loe$.

(ii) If $F$ is an open connected set,
$\partial F$ contains an arc $I$ of $\partial \E$, and one of the foci of
$\E$ is in
$F$, then $F$ is a strong peak set for
$L^1_h(\E)$.
}

\medskip
\noindent
{{\bf{Proof.}}}
(i) By [Sh1, p.21], there is a bounded function $w$ on $U$ such that
$$
\int_{\E} h dA \ = \int_{U} w h dA
$$
for all $h$ in $L^1_h(\E)$. So just as in the proof of the first half of
Theorem~5.10, we get
$$
|\int_F h dA | \ \leq \
(\| w \| + 1) \int_{\E \setminus F} |h |dA .
$$
(ii) If $F$ is not a strong peak set for $\loe$, as in Theorem~5.10
we can find a function $u \in C^1({\bf R}^2)$ that has $\Delta u = 1$ on $F$ and
vanishes along with its gradient on $I$. Therefore it coincides with the
modified Schwarz potential $v$ of $\partial \E$ on $F$.
But ${\partial v}\over{\partial z}$ has square root  type
branch points at $\pm 1$
[Sh1,p.21], so $v$ is not $C^1$ in any neighborhood of a focus.
\ep

\bigskip
Remark: The preceding theorem and proof remain valid for ellipsoids in
${\bf R}^n$, where the pair of foci are replaced by the
$(n-1)$-dimensional focal ellipsoid (or caustic). See [Kh1] and [Sh1].
\medskip

Interestingly, {\it any} neighborhood of a {rough} boundary point is
automatically a strong harmonic peak set.

\bigskip
\noindent
{{\bf{Theorem 5.12.}}} {\it
Let $G$ be a domain in $\rn$
and $F$ an open subset of $G$ such
that $\partial F \cap
\partial G \cap \partial  \overline{G}^c$
contains an $(n-1)$-dimensional manifold $J$.
If $F$ fails to be a strong $\loh$ peak set, then there is a function $u$,
in $C^{2 - \varepsilon} (\rn)$ for all $\varepsilon > 0$,
such that $u$ and $\nabla u$ vanish on $J$, but $u$ is not identically zero
in a neighborhood of any point on $J$. Moreover, if $n=2$, and $J$ is a
Jordan arc, then $J$ must
actually be an analytic arc.
}

\medskip
\noindent
{{\bf{Proof.}}}
As in the proof of Theorem 5.10, if $F$ is not a strong peak set, there
is a function $\psi$ in $L^\infty(G)$ that annihilates $\loh$ and equals
one on $F$.
Then $u =  E\ast(\psi) $
satisfies
equation (5.3).

In the case $n=2$, it follows from [Sh1,p.39] that the existence of $u$
satisfying (5.3) forces $J$ to be an analytic arc.
\ep

\bigskip
{\bf Example 5.13}: If $G$ is a square, it follows 
from Theorem 5.12 that any neighborhood of a
corner is a strong $\loh$ peak set.  
More is true: if $F$ is a ribbon connecting two different sides (though maybe
missing the corner), then it is still strong peak. This is because a $u$
satisfying equation (5.3) would actually be uniquely determined by knowing it
vanished along with its derivative on an arc of one side of the square - it
would have to be the modified Schwarz potential of a half-plane. But it would
also have to be the modified Schwarz potential of another half-plane,
corresponding to the other side that $F$ touches. These two functions are
different, and cannot agree on any open set.

However, if $F$ is a large set that only touches one side, it will not be a
strong peak set. For there is a $C^\i$ function $v$ on $\rn$, 
identically $1$ on a neighborhood of $F$, and identically zero on a
neighborhood of the three sides that $F$ doesn't touch. Let $u$ be the
modified Schwarz potential of the side $F$ does touch. Then it follows from
Green's theorem that $f = \Delta (uv)$ annihilates $\loh$; moreover 
$f$ is $1$ on $F$ and in $L^\i$, so $F$ can not be a strong peak set.

\medskip
Clearly the ideas in Example 5.13 could be extended to
other domains.

\bigskip\bigskip

\centerline{\big{6. A Proof of the AGHR Theorem}}
\bigskip
Our methods allow us to give  new proofs of the results of
Armitage, Gardiner,  Haussmann and Rogge [AGHR].

Let $\rho = \rho_n = 2^{-1/n}$, and let
$\BO$ be the open ball centered at zero of radius $\rho$ (so it has exactly
half
the volume of $\B$).
Let $\sigma $ be the function that is $-1$ on $\BO$, $+1$ on $\B \setminus
\BO$, and $0$ off $\B$.

For $n \geq 2$, let $\L$ be the differential operator
on $\R^n$ given by
$$
\L(f) \ = \ \sum_{j=1}^n x_j {{\partial f} \over{\partial x_j}} + 
{{n-2}\over{2}} f .
$$

First we prove the following Lemma.

\medskip
\noindent
{\bf Lemma 6.1.}
{\it
Suppose $g$ is in $\lib$ and $\| g \| \leq 1$. If $n=2$, suppose also that
$\int_\D g = 0$. Then for all $y$ in $\B$ with $| y | = \rho$, we
have
$$
| \L_y [E\ast g (y) ]| \leq |\L_y [E\ast \sigma(y)] | ,
$$
with strict inequality unless $g$ is, almost everywhere, a unimodular
constant times $\sigma$.
}

\medskip
\noindent
{\bf Proof.}
First assume $n \geq 3$. Then a calculation yields that
$$
\L_y | x- y|^{2-n} \ = \ \left( {{n-2}\over{2}} \right)
{{|x|^2 - | y|^2}\over{|x -y |^n}}.
$$
Therefore, as $E(x-y) = c | x- y|^{2-n} $ for the appropriate constant $c =
c(n)$,
we get that 
$$
\L_y \ E\ast g (y)  \ = \ {{n-2}\over{2}} \, c \int_B
{{|x|^2 - | y|^2}\over{|x -y |^n}} g(x) dx .
\eqno{(6.1)}
$$
For $| y | = \rho$, the right-hand side of (6.1) is maximized if 
$$
g(x) \ =\  sgn \left({{|x|^2 - | y|^2}\over{|x -y |^n}}\right)
\ = \ \sigma(x) .
$$
Moreover, there will be cancellation in the integral in (6.1) unless
$g$ is a unimodular constant times $\sigma$.

Now consider the case $n=2$.
A calculation gives
$$
\L_y [\log |y-x|^2] - 1 \ =\
{{ |y|^2 - |x|^2}\over{ |y-x|^2}} .
$$
Therefore
$$
\eqalignno{
\L_y \ E\ast g (y)  \ &=\
c \int_\D \L_y \log|y-x|^2 g(x) dx \cr
&=\ c \int_\D \left[ \L_y \log|y-x|^2 - 1 \right] g(x) dx \cr
&=\ c \int_\D \left[{{ |y|^2 - |x|^2}\over{ |y-x|^2}}
\right] g(x) dx 
\cr }
$$
As before, this will be maximized when $|y | = \rho$ 
by $g(x) = \sigma(x)$.
\ep
\bigskip

Now we can prove 
Proposition~2 from [AGHR].

\medskip
\noindent
{\bf Theorem 6.2.}
{\it
Let $F \subseteq \B$, and assume $F$ is open and connected, and
$\partial F$ contains a relatively open subset of $\partial \B$.
Suppose $g$ annihilates $\lob$, $\| g \|_{\i} = 1$
and $g \equiv 1$ on $F$.
Then $F$ has empty intersection with $\BO$.
}

\medskip
\noindent
{\bf Proof.}
Let $\O = F \cap (\B \setminus \overline{\BO})$. Then both $E\ast g$ and
$E\ast\sigma$ have 
Laplacian $1$ on $\O$, and vanish along with their gradients on
$\partial \O \cap \partial \B$ (since they both vanish identically off $\B$).
Therefore they agree on $\O$, and in particular $\L (E\ast g) 
= \L(E\ast\sigma)$ on $\O$. 
If $\partial \O \cap \partial \BO$ is non-empty, then Lemma 6.1 forces $ g$ to
equal $\sigma$. 
\ep
\bigskip

Note that in dimension $2$, one only needs $F \cap \partial \D$ to have
positive measure.

In the terminology of Section~5.2, Theorem~6.2 says that if $F$ is a domain
containing a full neighborhood of $\partial B$, then $B(F) \leq 1$ if and
only if $F \cap \B_0$ is empty.

We need the following result for the case that $G$ is the ball and $K$ the
center point. As we think it may be useful in other cases, we 
give it in greater generality. Note that hypothesis~(6.3) will be satisfied
if, for example, the capacity of $K$ is zero and $\mu$ is positive.

\medskip
\noindent
{\bf Proposition 6.3.}
{\it Suppose $G$ is a domain in $\rn$, with 
piecewise smooth boundary, that satisfies a quadrature identity
$$
\int_G h(x) dx \ = \ \int_K h(x) d\mu(x)
\eqno{(6.2)}
$$
for all $h$ in $L^1_h(G)$, where $\mu$ is a signed measure supported on
$K$, and $K$ 
has
$(n-1)$-dimensional Hausdorff measure zero. 
Let $\um = E \ast \mu$ be the Newtonian potential of $\mu$, and assume
$$
\limsup_{ \rn \setminus K\, \ni \, x \to y} [ |\um(x)| + |\nabla \um (x)| ] \
= \ \infty \qquad \forall y \in K.
\eqno{(6.3)}$$

Let $\o$ be continuous on $\overline{G}$ and subharmonic on $G$, and assume
that it is badly approximable in $L^1_h(G)$. 
Then if $\o$ is non-negative on $K$, it is non-negative on $G$.
}

\medskip
\noindent
{\bf Proof.}
As $\o$ is badly approximable, there is a function $g$ in the ball of
$L^\i(G)$ that agrees with $sgn(\o)$ when $\o \neq 0$ and that annihilates
$\loh$; let us extend this function to be $0$ off $G$, and denote the new
function also by $g$.
Let 
$$
\eqalignno{
u\  & = E \ast (\chi_G - \mu) \cr
v\  & = E \ast g
\cr}
$$
Then $v$ is $C^1$ on $\rn$, $u$ is $C^1$ on $\rn \setminus K$,
and both vanish identically off $G$.

Note first that if $P_0$ is any component of $ P := \{ \o > 0 \}$, 
then by subharmonicity and continuity 
of $\o$ we must have that $\partial P_0$ contains a
relatively open subset of $\partial G$. 
As $u$ and $v$ agree outside $G$, it follows from Holmgren's theorem
(which asserts that if a harmonic function and its gradient both vanish on an
$(n-1)$-dimensional manifold, then the function must be identically zero) and
the fact that $K$ has $(n-1)$-dimensional Hausdorff measure $0$ that the
function $u-v$, which is harmonic on
$P_0 \setminus K$, must vanish identically on $P_0 \setminus K$.
As $v$ and $u - \um$ are $C^1$, it follows from~(6.3) that 
$K$ must be disjoint from $\overline{P}$.

Now let $N$ be a component of $\{ \o < 0 \}$. By hypothesis, $K \cap N =
\emptyset$.
Moreover, as $\o$ is subharmonic, $\partial N \cap G \subseteq \partial P$.

Claim: $K$ is disjoint from $\overline{N}$.
 
(i) If $\partial N$ contains a relatively open subset of $\partial G$, then as
before the fact that $u+v$ is harmonic on $N$ and zero off $G$ forces it to be 
zero on $N$. Therefore (6.3) implies $K \cap \overline{N} = \emptyset$.

(ii) If $\partial N$ does not contain a relatively open subset of $\partial G$,
then $\partial N \cap G$ is dense in $\partial N$, so $\partial N \subseteq
\partial P$, and therefore $K$ is disjoint from $\overline{N}$.

Now consider $\nabla( u - v )$. This is a harmonic vector field on $N$,
continuous on $\overline{N}$. Moreover, it is zero on $\partial N$
(because it is zero on $\partial P$ and $\partial G$).
Therefore on $N$, the function $u-v$ is constant.
As $\Delta(u-v) =2$, this forces $N$ to be empty.
\ep

The main result of [AGHR] now follows from Theorem~(6.2) and Proposition~(6.3).

\medskip
\noindent
{\bf Corollary 6.4.}
{\it
Suppose $\o$ is continuous on $\overline{\B}$ and subharmonic on $\B$, and 
that $h$ is continuous on $\overline{\B}$ and harmonic on $\B$.
Then $h$ is a best $L^1$-approximant to $\o$ if and only if
\item{(i)} $h = \o$ on $\partial \BO$, and
\smallskip
\item{(ii)} $h \leq \o$ on $\overline{\B} \setminus \BO$.
}

\medskip
\noindent
{\bf Proof.} (Sufficiency) If hypotheses (i) and (ii) hold,
then $sgn(\o - h) = \sigma$ whenever $\o - h$ is non-zero.
As $\sigma$ annihilates $\lob$, it follows that $h$ is a best harmonic
approximant of $\o$.

(Necessity) Conversely, if $h$ is a best harmonic
approximant of $\o$, let $f = \o - h$.
As $f$ is badly approximable, there is a function $g$ of norm $1$ in 
$\lib$ that annihilates $\lob$ and agrees with $sgn(f)$
whenever $f$ is non-zero.

As $f$ is subharmonic and continuous on $\overline{\B}$, it cannot be
strictly positive at any point of $\BO$ without being positive on a set $F$
which satisfies the hypotheses of Theorem~6.2. So by that theorem, 
we must have that $f \leq 0$ on $\BO$.

By the sub-mean value property of subharmonic functions, we must also have
$$
\partial \{ f < 0 \} \cap \B \ 
\subseteq \ \partial \{ f > 0 \} .
$$

Therefore we must either have that $f < 0$ on $\BO$, or $f \equiv 0$ 
on $\BO$. In the first case, 
$g$ must equal $\sigma$ a.e., and (i) and (ii) follow.
In the second case, (i) is immediate, and (ii) follows from Proposition~6.3,
as $f(0) \geq 0$ forces $f$ to be non-negative on all of $\overline{B}$.
\ep

\medskip

Another consequence of Lemma~6.1 is the following ``equigravitational''
result, which was suggested to us by
Bj\"orn Gustafsson.

\medskip
\noindent
{\bf Corollary 6.5.}
{\it
Let $K \subseteq \overline{\B}$ be a closed set with volume equal to the
volume of $\BO$, and such that its potential $U_K := E\ast\chi_K$ agrees
outside $\B$ with $U_{\BO}$. If $K \neq \overline{\BO}$, then no boundary
point $y$ of $\BO$ can be joined to $\partial \B$ by an arc $\Gamma$ that is
disjoint from $K \setminus \{ y \}$.
}

\medskip
\noindent
{\bf Proof.} Define $g$ to be $-1$ on $K$ and $+1$ on 
$\B \setminus K$. 
If there were such an arc $\Gamma$, it could be thickened to give an open set
$F$ which does not meet $K$ except possibly at $y$.
As in the proof of Theorem~6.2, we have $E\ast g = E \ast \sigma$ in $F$,
and Lemma~6.1 gives a contradiction.
\ep

\bigskip
\bigskip
\centerline{\big{7. Smooth functions with unbounded best approximants }}
\bigskip
First we characterize the best harmonic approximant to the Newton kernel with
pole in the ball of radius $\rho_n^2$, where as before $\rho_n = 2^{-1/n}$.
For any point $y$ in $\rn$, let $y^\prime$ be the Kelvin reflection 
in the sphere $\partial \BO$, {\it i.e.} $y^\prime$ is on the same ray
through the origin as $y$ and $|y | | y^\prime | = \rho_n^2$.
We shall continue to use $\sigma$ to denote the function that is $-1$ on
$\B_0$ and $+1$ on $\B \setminus \B_0$.

\medskip
\noindent
{\bf Theorem 7.1.}
{\it For $n \geq 3$, the best harmonic approximant in $L^1(\B_n)$ of the
function $\displaystyle f(x) = {{1}\over{| x - y |^{n-2}}}$ when $|y | \leq
\rho_n^2$ is the function
$$
h(x)\  =\ \left( {{\rho_n} \over{ | y |}}\right)^{n-2} 
{{1}\over{| x - y^\prime |}^{n-2}} .
$$
For $n = 2$, the best $L^1_h(\D)$ approximant to 
$f(x) = \log | x - y |$ for $|y| \leq {{1}\over{2}}$ is the function
$$
h(x) \ =    \ \log \sqrt{2} |y||x-y^\prime|
$$
When $y = 0$, the best approximants are the constant functions $\displaystyle
{{1}\over{\rho_n}}$ and $\displaystyle \log {{1}\over{\sqrt{2}}}$
respectively.
}
\noindent
{\bf Proof.}
Let $|y| \leq \rho_n^2$.
By direct computation, 
$$
| x- y | \ < \ {{|y|}\over{\rho_n}} |x - y^\prime|
$$
if and only if $x$ is in $\BO$. So $sgn(f-h) = - \sigma$
and annihilates $\lob$, and therefore $h$ is the best harmonic approximant.
\ep
\bigskip
Notice that if $f$ is replaced by $\min(f, M)$ for some large constant $M$,
or even by a $C^\i$ smoothing, the function $sgn(h-f)$ will still be
$\sigma$, so $h$ will still be the best approximant (if $n=2$, take the
cut-off from below). Letting $|y| = \rho_n^2$, therefore, we get:

\medskip
\noindent
{\bf Corollary 7.2.}
{\it
There exists a 
$C^\i$-function that is real-analytic in a neighborhood of
$\partial \B$ and whose best harmonic approximant is unbounded on $\B$.
}
\bigskip

This is in marked contrast with the behaviour in $L^2$: 

\medskip
\noindent
{\bf Theorem 7.3.}
{\it
If $G$ is a domain in $\rn$ with smooth boundary that is real-analytic near
the boundary point $x_0$, and $f$ in $L^2(G)$ extends real-analytically
across $x_0$, then its best approximant in $L^2_h(G)$ also extends
real-analytically
across $x_0$. }

\noindent
{\bf Proof.}
Let $u$ be the orthogonal projection of $f$ onto $L^2_h(G)$, so 
$$
f \ =\ u+ g 
$$
where $g$ is in $L^2(G)$ and annihilates $L^2_h(G)$.
By the harmonic analogue of Khavin's Lemma, there is 
$v$ in $W_0^{2,2}(G)$ with $\Delta v = g$ in $G$. 

As $f$ extends real-analytically across $x_0$, and denoting the extension
also by $f$, there is, in some small ball $B$ centered at $x_0$, a solution
to the Cauchy problem
$$
\Delta w = f, \quad w = 0 = \nabla w \ {\rm on}\ \partial G \cap B .
$$
Let $\O = G \cap B$. Then on $\O$,
we have $\Delta w = u + \Delta v$, so $\Delta \Delta (w - v) = 0$.
Thus, $w-v$ satisfies the biharmonic equation in $\O$, and vanishes along
with its gradient on $\partial \O \cap \partial G$ ({\it i.e.} a trace of the
function in $W^{2,2}(\O)$ does). 

As $\partial \O \cap \partial G$ is is real-analytic near $x_0$, by
``regularity up to the boundary'' theorems for elliptic operators
[F, p.205] we get that $w-v$ extends real-analytically across $x_0$, and so
therefore does $v$. Thus we get that $u = f - \Delta v$ extends
real-analytically across $x_0$.
\ep

Another corollary to Theorem~7.1 is the following:

\noindent
{\bf Corollary 7.4.}
{\it
If $\|g\|_\i \leq 1$ and $g$ annihilates $\loh$, then
$$
| E \ast g (y) | \leq | E \ast \sigma (y) |, \quad |y| \leq
\rho_n^2 ,
$$
with strict inequality unless $g$ equals a.e. a unimodular constant
times
 $\sigma$. Moreover, $\rho_n^2$ is the largest number for which
this is
true. }

\medskip
\noindent
{\bf Proof.}
For simplicity, we give the proof in the case $n \geq 3$; the case
$n=2$ is
similar.
Let  $h_y$ be the best harmonic approximant to $\displaystyle
{{1}\over{|x-y|^{2-n}}}$. For $|y| \leq
\rho_n^2 $, we have 
$$
\eqalignno{
|E\ast g (y) | \ &= \ c \left | \int_\B  \left[ {{1}\over{|x-y|^{2-n}}} - h_y(x)
\right] g(x) dx \right| \cr
&\leq\ c \int _\B \left| {{1}\over{|x-y|^{2-n}}} - h_y(x)
 \right| dx \cr
&=\ c \int_\B \left[ h_y(x) - {{1}\over{|x-y|^{2-n}}} \right] \sigma(x) \cr
&=\ |E\ast \sigma(y)|
\cr}
$$
Clearly equality requires $g$ to be a constant times $\sigma$.

Now, if $|y| > \rho_n^2$, we cannot have
$$
sgn \left( [h_y(x) -{{1}\over{|x-y|^{2-n}}} ] \cdot \sigma(x) \right)
$$
constant a.e. For this would force 
$\displaystyle h_y(x) -{{1}\over{|x-y|^{2-n}}} $ to vanish on $\partial \BO$.
If $\rho_n^2 < |y| \leq \rho_n$, this would force $h_y$ to have a pole at
$y^\prime$ which is inside $\B$; and if $\rho_n < |y| < 1$, 
this would force $h_y$ to have a pole at $y$.

So if $\displaystyle s(x) = sgn( h_y(x) -{{1}\over{|x-y|^{2-n}}}) $,
then $|E\ast s(y) |$ will be strictly larger than 
$|E \ast \sigma(y)|$.
\ep

Similarly we have
\medskip
\noindent
{\bf Corollary 7.5.}
{\it Let $K \subseteq \overline{\B}$ be a closed set with volume equal to the
volume of $\BO$, and such that its potential $U_K := E\ast\chi_K$ agrees
outside $\B$ with $U_{\BO}$. If $K \neq \overline{\BO}$, then
$|U_K(y)| < |U_{\BO} (y)| $ for $|y| \leq \rho_n^2$.
}

\bigskip
Let us mention one last consequence of these ideas.
Let $y_0 \in \B$, thought
of as close to the boundary. Let $h(x)$ be the best harmonic approximant of
$E(y_0 -x)$ and $s(x)$ be $sgn[ E(y_0 -x) - h(x)]$.
Let $F$ be an open connected set such that $\partial F$ contains a relatively
open subset of $\partial B$, and with $y_0$ in $\overline{F}$.
Then if $g$ is in the closed unit ball of $\lib$, annihilates $\lob$
and equals $s$ on $F$, then $g$ must equal $s$ a.e. on $\B$.
For indeed, $E\ast g = E \ast s$ on $F$, so 
$$
E\ast g (y_0) \ = \ \int [E(y_0 - x) - h(x) ] g(x) \ = 
\ \int [E(y_0 - x) - h(x) ] s(x).
$$
Therefore there is no cancellation in the first integral, and so $g$ must equal
$s$ a.e.

In other words, knowledge of $g$ on the (small) set $F$, along with the fact
that $g$ annihilates $\lob$ and is of norm $1$, uniquely determines
it.

\bigskip
In the analytic case, we can construct a
continuous function with unbounded best approximant, but have not been able
to make $\o$ any smoother:

\noindent
{\bf Proposition 7.6.}
{\it There is a function $\o$ that is continuous on the closed disk and whose
best
analytic approximant in $L^1(\D)$ is unbounded near every point of $\partial
\D$.}

\medskip
\noindent
{\bf Proof.}
Let $f = u + iv$ be a holomorphic function on the unit disk, whose imaginary part is
continuous on $\partial \D$ and whose real part is positive and unbounded on
$\partial \D$ ({\it e.g.} the Riemann map onto the set $\{ x + i y : x > 1,
0 < y < {{1}\over{x}} \}$).
By taking a suitable convex combination of rotates of $f$, we can moreover
assume that $u$ is unbounded near every point of $\partial \D$, and that $f$
is in $A^1$.

Let 
$$
\o(z) \ = \ 2(1 - |z|^2) u(z) + i v(z) .
$$
Then $\o$ is continuous on $\overline{\D}$, because $u(z) = o(\log |1-z|)$.
Moreover,
$$
\o(z) - f(z) \ = \ [1 - 2 |z|^2] u(z) ,
$$
which is positive in $\D_0$ and negative outside $\D_0$.
Therefore $f$ is the best analytic aproximant to $\o$.
\ep

\bigskip
{\bf Question} If $\omega$ is H\"older continuous on $\overline{\D}$ 
must its best $A^1$ approximant be continuous on $\overline{\D}$?


\bigskip\bigskip

\centerline{{\big{References}}} \par

\bigskip
\item{{\bf{[Ak] :}}}     N.I. Akhieser, {{\it{Lectures in Approximation
                         Theory}}}, Nauka, \hbox{Moscow}, 1965 \hfill \break
                         \hbox{(in Russian)}.
\smallskip
\item{{\bf{[AGHR] :}}}   D.H. Armitage, S.J. Gardiner, W. Haussmann, and
			 L. Rogge, Characterization of best harmonic and
			 superharmonic $L^1$-approximants, {\it J. Reine
			 Angew. Math.} {\bf 478}(1996), 1-15.
\smallskip
\item{{\bf{[BW] :}}}     J. Bourgain and T. Wolff, Note on gradients of
			 harmonic functions in dimension $\geq 3$, 
			 {\it Colloq. Math.} {\bf 60/61} 1 (1990), 253-260.

\smallskip
\item{{\bf{[CJ] :}}}     L. Carleson and S. Jacobs, Best uniform
                         approximation by analytic functions, {{\it{Ark.
                         Mat.}}} {{{\bf{10}}}(1972), \hbox{219--229}.
\smallskip
\item{{\bf{[D] :}}}      P. Duren, {{\it{The Theory of
                         \hbox{$H^{p}$-Spaces}}}}, Academic Press,
                         \hbox{New York}, 1980.
\smallskip
\item{{\bf{[DKSS] :}}}   P. Duren, D. Khavinson, H.S. Shapiro, and
                         C. Sundberg, Contractive zero divisors in
                         \hbox{Bergman} spaces, {{\it{Pacific J. Math.}}},
                         {{\bf{157}}}(1993), \hbox{37--56}.
\smallskip
\item{{\bf{[EG] :}}}      L.C. Evans and R.F. Gariepy, {{\it{Measure theory 
			 and fine properties of functions}},
			 CRC Press, \hbox{Boca Raton}, 1992.
\smallskip
\item{{\bf{[F] :}}}      A. Friedman {\it Partial Differential Equations},
			 Holt, Rinehrat and Winston, New York, 1969.

\smallskip
\item{{\bf{[GK] :}}}     T. Gamelin and  D. Khavinson, 
			 The isoperimetric inequality and rational
			 approximation, {\it Amer. Math. Monthly},
			 {\bf 96} (1989) No. 1, \hbox{18--30}.
\smallskip
\item{{\bf{[GHJ] :}}}    M. Goldstein, W. Haussmann, and K. Jetter, Best
                         harmonic \hbox{$L^{1}$-approximation} to
                         subharmonic functions, {{\it{J. \hbox{London} Math.
                         Soc.}}} (2), {{\bf{30}}}(1984), \hbox{257--264}.
\smallskip
\item{{\bf{[GT] :}}}     D. Gilbarg and N.S. Trudinger, {\it Elliptic Partial
			 Differential Equations of Second Order,}
			 Springer, Berlin, 1977.
\smallskip
\item{{\bf{[HKL] :}}}    W.K. Hayman, D. Kershaw, and T.J. Lyons, The best
                         harmonic approximation to a continuous function,
                         Proceedings of the Conference on Functional
                         Analysis and Approximation, \hbox{Oberwolfach}
                         1983, {{\it{International Series on Numerical
                         Mathematics}}} {{\bf{65}}}(1984), (eds.
                         \hbox{P.L. Butzer}, \hbox{B.Sz.-Nagy}, and
                         \hbox{R. Stens}) \hbox{Birkh{\"{a}}user},
                         \hbox{Basel}, \hbox{317--327}.
\smallskip
\item{{\bf{[HS] :}}}     E. Hewitt and K. Stromberg, {{\it{Real and Abstract
                         Analysis}}}, \hbox{Springer}--\hbox{Verlag}, \hfill
                         \break \hbox{New York}, 1965.
\smallskip
\item{{\bf{[Ka] :}}}     J.-P. Kahane, Best approximation in $L^{1} \left(
                         {{\bf{T}}} \right) $, {{\it{Bull. Amer. Math.
                         Soc.}}}, \hfill \break {{\bf{80}}}(1974),
                         \hbox{788--804}.
\smallskip
\item{{\bf{[Kh1] :}}}    D. Khavinson {\it Holomorphic Partial
			 Differential Equations and Classical Potential Theory,}
                         \hbox{Universidad de la Laguna}, 1996.
\smallskip
\item{{\bf{[Kh2] :}}}    S.Ya. Khavinson, Foundations of the theory of
                         extremal problems for bounded analytic functions
                         and various generalizations of them, {{\it{Amer.
                         Math. Soc. Transl.}}} (2),
                         \hbox{Vol.{{\bf{129}}}(1986)}, \hbox{1--56}.
\smallskip
\item{{\bf{[Kh3] :}}}    S.Ya. Khavinson, Foundations of the theory of
                         extremal problems for bounded analytic functions
                         with additional conditions, {{\it{Amer. Math. Soc.
                         Transl.}}} (2), Vol. {\bf{129}}(1986),
                         \hbox{64--114}.
\smallskip
\item{{\bf{[Kh4] :}}}    S.Ya. Khavinson, Duality relations and criteria
                         for elements giving the best approximation,
                         {{\it{Lecture Notes}}}, \hbox{Moscow} Institute of
                         Civil Engineering, \hbox{Moscow}, 1976,
                         \hbox{1--46} \hbox{(in Russian)}.
\smallskip
\item{{\bf{[Kh5] :}}}    S.Ya. Khavinson, \hbox{Chebyshev's} ideas in the
                         theory of best \hbox{approximation, II, III},
                         {{\it{Lecture Notes}}}, \hbox{Moscow} Institute of
                         Civil Engineering, \hbox{Moscow}, 1977,
                         \hbox{1--20} and \hbox{1--21} \hbox{(in Russian)},
                         (\hbox{ed. by} \hbox{E.Sh. Chatskaya}).
\item{}                  (\hbox{Part II:}   Completeness of Systems.)
\item{}                  (\hbox{Part III:}  Uniqueness and non-uniqueness of
                         best approximations.)
\smallskip
\item{{\bf{[Kh6] :}}}    S.Ya. Khavinson, On uniqueness of the function
                         giving the best approximation in the metric of the
                         space $L^{1}$, {{\it{Izv. Akad. Nauk SSSR}}}, ser.
                         matem. {{\bf{22}}}(1958), \hbox{243--270}
                         \hbox{(in Russian)}.
\smallskip
\item{{\bf{[KP-GS] :}}}  D. Khavinson, F. Per{\'{e}}z-Gonz{\'{a}}lez, and
                         H.S. Shapiro, Approximation in \hfill \break
                         \hbox{$L^{1}$-norm} by elements of a uniform
                         algebra, {\it Constructive Approximation,} to
                         appear
\smallskip
\item{{\bf{[KS] :}}}     D. Khavinson and M. Stessin, Certain linear
                         extremal problems in \hbox{Bergman} spaces of
                         analytic functions, {\it Indiana Math. J.,} 
			 {\bf 46} (1997), \hbox{933--973}
\smallskip
\item{{\bf{[RS] :}}}     W.W. Rogosinski and H.S. Shapiro, On certain
                         extremum problems for analytic functions,
                         {{\it{Acta Math.}}} {{\bf{90}}}(1953),
                         \hbox{287--318}.
\smallskip
\item{{\bf{[R] :}}}      V.G. Ryabych, Extremal problems for summable
                         analytic functions, {{\it{Siberian Math. J.}}}
                         {{\bf{XXVIII}}}(1986), \hbox{No.3},
                         \hbox{212--217}, \hbox{(in Russian)}.
\smallskip
\item{{\bf{[Sh1] :}}}    H.S. Shapiro, {{\it{The \hbox{Schwarz} Function and
                         its Generalizations to Higher Dimensions}}},
                         \hbox{Wiley}, \hbox{New York}, 1991.
\smallskip
\item{{\bf{[Sh2] :}}}    H.S. Shapiro, Regularity properties of the element
                         of closest approximation, {{\it{Trans. Amer. Math.
                         Soc.}}}, {{\bf{181}}}(1973), \hbox{127--142}.
\smallskip
\item{{\bf{[V] :}}}      D. Vukoti\'c, Linear extremal problems for
   			 Bergman spaces, {\it Exposition. Math.}, 
			 {\bf 14} (1996) No. 4, 313--352.
\smallskip
\item{{\bf{[W] :}}}      J. Walsh, {{\it{Interpolation and Approximation by
                         Rational Functions in the Complex Domain}}},
                         \hbox{American} Mathematical Society,
                         \hbox{Providence, Rhode Island}, 1960.

\bigskip
\noindent
D. Khavinson \hfill \break
Department of Mathematics \hfill \break
University of Arkansas \hfill \break
Fayetteville, Arkansas~~72701 \hfill \break
U.S.A. \hfill \break
e-mail:  dmitry@comp.uark.edu \par

\bigskip
\noindent
J.E. M\raise.45ex\hbox{c}Carthy \hfill \break
Department of Mathematics \hfill \break
Washington University \hfill \break
Saint Louis, Missouri~~63130 \hfill \break
U.S.A. \hfill \break
e-mail:  mccarthy@math.wustl.edu \par

\bigskip
\noindent
H.S. Shapiro \hfill \break
Department of Mathematics \hfill \break
Royal Institute of Technology \hfill \break
Stockholm,~~S-10044 \hfill \break
SWEDEN \hfill \break
e-mail:  shapiro@math.kth.se 

\end